\documentclass[10pt, final]{amsart} 
\usepackage[english]{babel}
\usepackage{dsfont}
\usepackage{times}
\usepackage{enumerate}
\usepackage{graphicx, graphics, float}
\usepackage{color}
\usepackage[usenames,dvipsnames]{xcolor}
\usepackage[colorlinks=true,urlcolor=blue,
citecolor=ForestGreen,linkcolor=blue,linktocpage,pdfpagelabels,
bookmarksnumbered,bookmarksopen,backref]{hyperref}
\usepackage{amssymb}
\usepackage{mathabx} 

\numberwithin{equation}{section}

\newtheorem{theorem}{Theorem}[section]
\newtheorem*{acknowledgement}{Acknowledgement}

\newtheorem{corollary}[theorem]{Corollary}

\newtheorem{lemma}[theorem]{Lemma}

\newtheorem{proposition}[theorem]{Proposition}

\theoremstyle{remark}
\newtheorem{definition}[theorem]{Definition}

\theoremstyle{remark}
\swapnumbers
\newtheorem{remark}[theorem]{Remark}

\newcommand{\R}{\mathds R} 
\DeclareMathOperator{\cat}{cat}

\DeclareMathOperator{\Crit}{Crit}
\DeclareMathOperator{\diam}{diam}
\DeclareMathOperator{\Vol}{Vol}
\DeclareMathOperator{\inj}{inj}
\DeclareMathOperator{\dive}{div}
\DeclareMathOperator{\vol}{vol}

\title[Van der Waals-Cahn-Hilliard equation with volume constraint]{Lusternik-Schnirelman and Morse Theory for the\\ Van der Waals-Cahn-Hilliard equation\\ with volume constraint}
\author[V. Benci, S. Nardulli, L. E. Osorio, P. Piccione]{Vieri Benci, Stefano Nardulli,\\ Luis Eduardo Osorio Acevedo\and Paolo Piccione}
\address{\begin{tabular}{lll}
Universit\`a di Pisa & & Universidade Federal do ABC\\
Dipartimento di Matematica & & Centro de Matem\'atica Cogni\c{c}\~ao Computa\c{c}\~ao\\
Via Filippo Buonarroti 1/c  & & Avenida dos Estados, 5001\\
56123 -- Pisa & &  Santo Andr\'e, SP, CEP 09210-580\\
Italy & & Brazil\\
\emph{E-mail}: {\tt vieri.benci@unipi.it} & &\emph{E-mail}: {\tt stefano.nardulli@ufabc.edu.br}
\\[.5cm]  Universidad Nacional de Colombia & & Universidade de S\~ao Paulo\\ Departamento de Matem\'atica y Estad\'istica & & Departamento de Matem\'atica\\ 
 Cra 27 No. 64-60, 170003 & & Rua do Mat\~ao 1010\\
Manizales, Colombia & & S\~ao Paulo, SP 05508--090, Brazil\\ \emph{E-mail}: {\tt leosorioa@unal.edu.co} & & \emph{E-mail}: {\tt paolo.piccione@usp.br}
\end{tabular}
}
\date{July 13, 2020}
\begin{document}
\begin{abstract}
We give a multiplicity result for solutions of the Van der Waals-Cahn-Hilliard two phase transition equation with volume constraints on a closed Riemannian manifold. Our proof employs some results from the classical Lusternik--Schnirelman and Morse theory, together with a technique, the so-called \emph{photography method}, which allows us to obtain lower bounds on the number of solutions in terms of topological invariants of the underlying manifold. The setup for the photography method employs recent results from Riemannian isoperimetry for small volumes.
\end{abstract}

\maketitle

\renewcommand{\contentsline}[4]{\csname nuova#1\endcsname{#2}{#3}{#4}}
\newcommand{\nuovasection}[3]{\hbox to \hsize{\vbox{\advance\hsize by -1cm\baselineskip=10pt\parfillskip=0pt\leftskip=3.5cm\noindent\hskip -2.5cm {#1}\leaders\hbox{.}\hfil\hfil\par}$\,$#2\hfil}}
\newcommand{\nuovasubsection}[3]{\hbox to \hsize{\vbox{\advance\hsize by -1cm\baselineskip=10pt\parfillskip=0pt\leftskip=3.9cm\noindent\hskip -2cm \textsl{#1}\leaders\hbox{.}\hfil\hfil\par}$\,$#2\hfil}}
\newcommand{\nuovasubsubsection}[3]{\hbox to \hsize{\vbox{\advance\hsize by -1cm\baselineskip=10pt\parfillskip=0pt\leftskip=4cm\noindent
\hskip -2cm #1\leaders\hbox{.}\hfil\hfil\par}$\,$#2\hfil}}

\tableofcontents

\section{Introduction}
The Van der Waals-Cahn-Hilliard two phase transition equation \eqref{Eq:EulerLagrangeVEpsilon} has attracted the interest of physicists, analysts, and geometers. 
This is a variational equation, obtained as the Euler-Lagrange equation of the energy functional $E_{\varepsilon}$ defined in equation \eqref{Eq:DefinitionOfFunctional}, which has
the classical form of a kinetic term plus a double well potential. 
Historically, the energy $E_{\varepsilon}$ was already proposed in 1898 by Van der Waals in \cite{vanderWaals} for the transition liquid-vapor phase. In 1958, 
J. W. Cahn and J.E. Hilliard in \cite{CahnHilliard1958} used it to model the transition of phase in some binary alloy.  Ginzburg-Landau used the same functional
to model ferromagnetic behaviour of materials \cite{Presutti2009}. 
A quick search through the literature shows the ubiquitous nature of this equation. For instance, for applications to Biology the interested reader can consult the book of Murray \cite{MurrayBook}. 

From the theoretical physics viewpoint, various attempts were made to derive equation \eqref{Eq:EulerLagrangeVEpsilon} as limit of some microscopically consistent model of statistical mechanics, however, this goal is yet to be completely achieved. 
Along this line,  in the papers \cite{GiacominLebowitz1}, \cite{GiacominLebowitz2}  a more complicated non-local equation is derived from a microscopical model whose properties have
analogies to the properties of  equation \eqref{Eq:EulerLagrangeVEpsilon}, and which constitutes a first order approximation.  Let us also mention the paper \cite{BertiniLandimOlla}, in which \eqref{Eq:EulerLagrangeVEpsilon} is obtained as the hydrodynamic limit of the Ginzburg Landau equation. For a complete account of statistical mechanics studies of these problems the interested reader can also consider the book \cite{Presutti2009}. 

A quite complete mathematical analysis of the positive functions realizing the minimum of $E_{\varepsilon}$ was carried out in the work of Modica \cite{Modica}. Following this seminal paper, many other authors gave results about the minimization of $E_{\varepsilon}$, see for example \cite{Baldo} for the case of multicomponent mixtures etc. 

In the present paper we consider a constrained variational problem for the functional $E_\varepsilon$. More precisely, 
we develop techniques to determine a family $u_{\varepsilon}$, $\varepsilon\in\left]0,\varepsilon_0\right]$, of critical points of $E_\varepsilon$ under the volume constraint $\int_Mu_{\varepsilon}=V$, with $V>0$ fixed and small, i.e., each $u_\varepsilon$ is a solution of problem \eqref{Eq:EulerLagrangeVEpsilon} below.
It follows from \cite[Theorem~1]{HutchinsonTonegawa} that, under mild conditions, $u_\varepsilon$ converges as $\varepsilon$ goes to $0$, in an appropriate geometric measure theoretic  sense, to a characteristic function of a finite perimeter set. This set has reduced boundary whose regular part  is a constant mean curvature (CMC) smooth hypersurface, and it is relatively open and dense into this boundary. However, in this singular limit, multiplicity issues may occur, and this affects the optimal regularity of the limit. 
It was only recently that a more advanced regularity theory became available, thanks to the works \cite{BW2018Stable} and \cite{BCW19}. To our knowledge, this more advanced regularity theory is only available for the unconstrained problem (as in \cite{ChodoshMantoulidis2018Minimal}). In this setting, it provides an alternative to the regularity of the limit interface, after the work of Y. Tonegawa and N. Wickramasekera, and based on Wickramasekera's regularity theory for stable minimal hypersurfaces, which was generalized to stable CMC hypersurfaces in \cite{BW2018Stable}, \cite{BCW19}. 
Under a different perspective, in \cite{PacardRitore} Pacard and Ritor\'e showed that every smooth constant mean curvature boundary can be suitably approximated by solutions of \eqref{Eq:EulerLagrangeVEpsilon}. 
These circumstances open new perspectives as to the search of CMC boundaries using this PDE approach, which is the objective pursued in this research project. We will address the limit procedure and the geometric consequences in a forthcoming paper. As a partial result in this direction it is worth mentioning that using the very recent work \cite{BW20}, it is immediate to show that a sequence of critical points with low-energy produced by our method converges to at least an almost embedded CMC boundary in the sense of \cite{BW20} enclosing a small volume. This is possible because our solution at low-energy have uniformly bounded energy and have uniformly bounded from above Morse index, by the dimension of the ambient manifold fulfilling the conditions required by \cite{BW20}. In this way we give an alternative proof of existence of at least one almost embedded CMC boundaries in the small volume regime. It seems to us that the theory of convergence developed in \cite{BW20} does not permit easily to maintain the lower bound on the solutions when we pass to the limit to get the same lower bound on the limit CMC almost embedded boundaries.

In view of Yau's conjecture for minimal hypersurfaces and the result of Pacard-Ritoré, it is reasonable to expect that our lower bound estimate is not optimal. It is conceivable that minimax methods may prove the existence of infinitely many solutions of the Cahn-Hilliard equation for small values of the temperature parameter $\varepsilon$ and any values of volume, using the limit problem for minimal (or CMC) hypersurfaces. However, this approach would require proof of Yau's conjecture for CMC boundaries, which at present is unforeseen, and also a further assumption on nondegeneracy, which only holds in a generic case. For a minimal (or CMC) hypersurface, nondegeneracy means that there are no nontrivial Jacobi fields in the kernel of the linearized operator associated with the second variation of the area functional. Such assumption is needed to apply the result of Pacard and Ritoré, which provides the link backward from a minimal (CMC) hypersurface and solutions of the Cahn-Hilliard equation.
For what concerns related results of existence of CMC boundaries by other methods it is worth to mention \cite{PacardXu, NarAnn, NarCalcVar} using the implicit function theorem and the isoperimetric problem for small volumes, \cite{arXiv:1808.03527v1,arXiv:1910.00989,arXiv:2004.05120,MR4011704} using the min-max construction. 
	The convergence of the Allen--Cahn--Hilliard energy to the perimeter functional was recently used to construct minimal hypersurfaces in any closed Riemannian manifold \cite{MR3743704,MR3814054} as an alternative approach to min-max methods  \cite{MR626027}. 
    The main results of the present paper are built upon a theory of multiplicity of solutions for semi-linear variational elliptic equations based on topological and nonlinear methods, along the lines
of \cite{BenciCerami}, \cite{BenciCeramiPassaseo}, \cite{BenciCeramiCalcVar}, \cite{BenciNewApproach}, \cite{BenciBonannoMicheletti}. In the last series of papers the authors treated a quite different kind of nonlinear problems.  A first essay to apply to the Cahn-Hilliard equation these abstract topological argument was achieved in a previous work (of three among the four authors of the present work) in \cite{BeNaPiCalcVar2020}. However in \cite{BeNaPiCalcVar2020} was considered a different setting, i.e., the potential had absolute minima at different levels (what we called asymmetric potential) the problem was studied in a compact set of the Euclidean space and the way in which the hypothesis of the abstract method were proved were completely different. In fact they were involving the theory and regularity of some suited auxiliary variational inequalities not related to the methods employed in the present paper. In this paper we treat for the first time in the literature (at the best of our knowledge) the Cahn-Hilliard equation on a compact manifold (having symmetric potential) with the afore mentioned topological methods combining with a fine understanding of the Riemannian isoperimetric problem for small volumes which use tools from geometric measure theory. More precisely, in order to establish a lower bound on the number of solutions of problem \eqref{Eq:EulerLagrangeVEpsilon} we employ a method from Lyusternik--Schnirelman and Morse theory, that will
be referred to as the \emph{photography method}, see Section~\ref{sec:photography} for details. Roughly speaking, a lower bound on the number of solutions that belong to a suitable
sublevel of the associated energy is given in terms of topological invariants of the underlying manifold. A correspondence between the topological invariants of the energy sublevel and
those of the underlying (finite-dimensional) manifold is produced by two continuous maps going in both ways, and whose composition is a homotopy equivalence of the finite-dimensional
manifold. The map from the finite-dimensional manifold to the sublevel is a sort of \emph{photography map}, which associates to each point a bell-shaped function around the point. This map reproduces a copy (the \emph{photography}) of the underlying manifold inside the energy sublevel. The map going backwards, i.e., from the sublevel to the finite-dimensional manifold, is given by a \emph{barycenter} map, which associates to each function, a suitably defined point in the domain around which most of the mass of the function is concentrated.
This construction is interesting when it can be made in such a way that the barycenter of a photography map is the identity map, up to homotopies. In this case, by an elementary topological argument the
Lysternik--Schnirelman category and each Betti number takes on the energy sublevel a larger value than the value it takes on the domain manifold, and the desired estimate follows from standard variational theories. We use the photography method to prove the existence of multiple solutions of our constrained variational problem; such solutions come in two classes: low energy solutions, and high energy solutions. By high energy solutions, we mean that we do not have an \emph{a priori} estimate of the energy. The result is obtained under suitable assumptions on the potential function $W$ associated to the problem, including a certain growth condition at infinity, see Section~\ref{sec:formulation} for details. It is important to observe that such asymptotic growth assumption can be weakened significantly in order to obtain the existence of low energy solutions. This will be discussed in Section~\ref{sec:dropping}.
\smallskip

The paper is organized as follows. Section~\ref{sec:formulation} contains the formulation of the PDE problem, with all technical assumptions on the potential function needed for the variational setup, and the statement of our main results. Section~\ref{sec:photography} is the core of the paper. After recalling some generalities on Lusternik--Schnirelman theory and Morse theory, we give a detailed description of the photography method, and its concrete application for the variational problem considered here. We prove the Palais--Smale condition, and we establish the properties of the photography map and the barycenter map. As to the photography map, our definition relies heavily on 
some geometric measure theoretical result proven by Modica in \cite{Modica} in the case of domains of $\mathds R^n$. For the development of our theory, we will need a formulation of the results in the context of Riemannian manifolds, and the details of this
formulation are given in Section~\ref{sec:geommetricpreliminaries} and in Proposition~\ref{Prop:1and3Modica}.
We point out that some formulations of Modica's result in the context of Riemannian manifolds have previously appeared in the literature - for the nonstationary case, see for instance Proposition $5.10$ in \cite{Pisante2016}. For the barycenter map, we employ a non-intrinsic approach by resorting to Nash embedding theorem, and we use several
extrinsic Riemannian geometry results combined with results obtained by the second author in a series of papers culminated in \cite{NardulliOsorioIMRN} co-authored by the forth author too. In particular, our approach requires
several technical results from isoperimetric theory that establish an estimate on the diameter of isoperimetric regions of small volume as can be found in \cite{MorganJohnson2000} and in \cite{NardulliOsorioIMRN}.   
\section{Formulation of the problem and main results}\label{sec:formulation}
In this section we will give the description of the nonlinear PDE problem, and we will formulate the main result concerning the multiplicity of its solutions. Let us assume that $W\colon\mathds R\to\left[0,+\infty\right[$ is a function of class $C^2$ and that $(M,g)$ is an $N$-dimensional compact Riemannian manifold without boundary; precise assumptions on $W$ will be given below. For fixed $\varepsilon,V>0$, we are concerned with the existence of multiple pairs $(u_{\varepsilon,V}, \lambda_{\varepsilon,V})\in H^1(M)\times\R$  such that the following equalities are satisfied:
\begin{equation}\label{Eq:EulerLagrangeVEpsilon}
\begin{aligned}
-\varepsilon\,\Delta u_{\varepsilon,V}+\tfrac{1}{\varepsilon}W^{\prime }(u_{\varepsilon,V}) &=\lambda_{\varepsilon,V}, \\
\int_{M}u_{\varepsilon,V}\,\mathrm dv_g &=V.
\end{aligned}
\end{equation}
As to the assumptions on $W$, we will consider a double well potential, i.e., a map satisfying the following assumptions: 
\begin{itemize}
\item[\textbf{(a)}] $W(s)$ has only two global minima, at $s=0$ and at $s=1$, and a unique local maximum at $s=\frac12$; moreover
\begin{equation}\label{Eq:Potential}
 W(0)=W'(0)=W(1)=W'(1)=0;\quad W''(0), W''(1)>0;
 \end{equation}
\item[\textbf{(b)}] there exist positive constants $A,B$ such that 
\begin{equation}\label{Eq:Potential0}
\big|W'(s)\big|\le A+Bs^{p-1},\quad\text{for some} \ p<\frac{2N}{N-2}=:2^*,\quad (p<\infty \ \text{if}\ \ N=1,2);
\end{equation}
\item[\textbf{(c)}] for some $\delta>0$:
\begin{equation}\label{Eq:WeakUpperBarriers}
W'(s)>0,\quad \forall\, s\in\left]1, 1+\delta\right];
\end{equation}
\item[\textbf{(d)}] there exists $c_1,c_2, t_0>0$ such that:
\begin{equation}\label{Eq:aPrioriL1Convergence+}
c_1|t|^{p_1}<W(t)<c_2|t|^{p_2},\quad \text{ when } |t|\ge t_0,
\end{equation}
and where $2<p_1<\hat2$, $p_1\le p_2\le 2(p_1-1)$, with $\hat2=\frac12\,{2^*}+1\le2^*+1$. 
\end{itemize}
\begin{figure}[H]
\begin{center}
\includegraphics[]{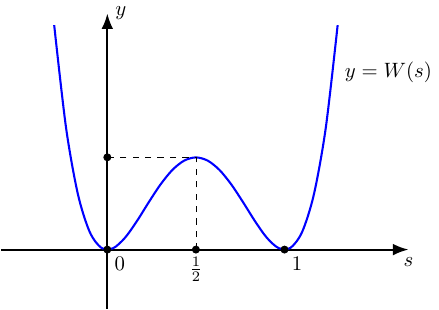}
\caption{The typical shape of the symmetric \emph{double well} potential $W$ considered in Problem \eqref{Eq:EulerLagrangeVEpsilon}.}
\end{center}
\end{figure}
The solutions of Problem \eqref{Eq:EulerLagrangeVEpsilon} are the critical points of the following energy functional
\[E_\varepsilon:H^1(M)\longrightarrow\R,\]
\begin{equation}\label{Eq:DefinitionOfFunctional}
E_{\varepsilon}(u)=\frac\varepsilon2\int_{M} \left\vert
\nabla u\right\vert^2\,\mathrm dv_g+\frac{1}{\varepsilon}\int_{M} W\big(u(x)\big)\,\mathrm dv_g,
\end{equation}
under the constraint 
\begin{equation*}
\int_{M} u\,\mathrm dv_g=V.
\end{equation*}
Here, $\mathrm dv_g$ denotes the volume density of the metric $g$.\medskip

Consider the following topological invariants of the manifold $M$.
Given a topological space $\mathcal X$, let us recall the definition of some topological invariants of $\mathcal X$:
\begin{itemize}
          \item $\cat(\mathcal X)$ is the Lusternik--Schnirelman category of $\mathcal X$, see Definition \ref{Def:Lusternik-SchnirelmannCategory},
          \item $\beta_k(\mathcal X)$ is the $k$-th Betti number\footnote{%
          Recall that the $k$-th Betti number of $\mathcal X$  is the dimension of the  $k$-th Alexander-Spanier cohomology vector space of $\mathcal X$ with coefficients in $\R$} of $\mathcal X$. Similarly, if $\mathcal Y\subset\mathcal X$ is a subspace,
          $\beta_k(\mathcal X,\mathcal Y)$ is the $k$-th Betti number of the pair;
          \item $P_1(\mathcal X)=\sum_k\beta_k(\mathcal X)$; this is the value at $1$ of the Poincar\'e polynomial of $\mathcal X$ (see Definition~\ref{thm:defpoincarepolynomial}).
          \end{itemize}
The main result of the paper gives a lower bound on the number of solutions of Problem \eqref{Eq:EulerLagrangeVEpsilon} in terms of these topological invariants of $M$.
\begin{theorem}\label{Thm:Main1}
Let $W$ satisfy assumptions (a), (b), (c) and (d) above. Then, there exists $V^*=V^*(M,g)>0$ such that for every $V\in\left]0,V^*\right[$ there exists $\varepsilon^*=\varepsilon^*(V, M,g, W)>0$, such that for every $\varepsilon\in\left]0,\varepsilon^*\right[$, Problem \eqref{Eq:EulerLagrangeVEpsilon} admits
at least $\cat(M)+1$ distinct solutions.
Moreover, if for some given $V$ and $\varepsilon$ as above all the solutions of Problem \eqref{Eq:EulerLagrangeVEpsilon} are \emph{nondegenerate} (i.e., they correspond to nondegenerate critical points of $E_\varepsilon$) then there are  at least $2P_{1}(M)-1$ solutions.
\end{theorem}
The nondegeneracy assumption in the last part of the statement can be omitted, provided that a suitable notion of multiplicity of solutions is taken into consideration, see  Definition~\ref{Eq:DefTopologicallyNondegenerate}.
\subsection*{Notations and terminology}
Given a Riemannian manifold $(M,g)$, we will denote by $N\in\mathbb{N}\setminus\{0\}$ the dimension of $M$, $\vol_g$ the volume function of the metric $g$, by $\inj_M$ the injectivity radius of $(M,g)$ (which is well defined and strictly positive when $M$ is compact), $\mathcal{H}_{g}^{N-1}$ is the $(N-1)$-dimensional Hausdorff measure associated to the distance induced on $M$ by the metric $g$,  and by $\diam_g$ the diameter function of sets induced by the metric associated to $g$.
\begin{acknowledgement} The authors wish to thank the anonymous referees for their careful proofreading of the original manuscript. This helped a lot in improving the quality of the presentation of our results. The second author is partially sponsored by "Bolsa no Exterior-Pesquisa" (Fapesp BPE, 2018/22938-4), "Jovens Pesquisadores em Centros Emergentes" (JP-FAPESP, 21/05256-0), and by "Bolsa de Produtividade em Pesquisa 2" (CNPq, 302717/2017-0), Brazil. The third author is sponsored by Fapesp (Postdoc Scholarship 2017/13155-3) and Minciencias Colombia (Postdoctoral Project 80740-738-2019). The fourth author is partially sponsored by Fapesp (Thematic Project, 2016/23746-6), and by "Bolsa de Produtividade em Pesquisa 1A" (CNPq, 313773/2021-1), Brazil.
The authors wish to thank João Henrique Santos de Andrade for the proof-reading of the final manuscript. 
\end{acknowledgement}
\section{Geometric measure theoretical preliminaries}
\label{sec:geommetricpreliminaries}
For the development of our theory, we will need a Riemannian counterpart of some results that were originally obtained by 
Modica in \cite{Modica} in the case of domains in Euclidean spaces. Although Modica's main ideas carry over to the geometrical setup
without major difficulties, for the reader's convenience we will give here a detailed proof of \cite[Proposition~2]{Modica} stated in our Riemannian context. 
Let us first recall some definitions.
\begin{definition} Let $(M,g)$ be a Riemannian manifold of dimension $N$, $U\subseteq M$ an open subset, $\mathfrak{X}_c(U)$ the set of smooth vector fields with compact support on $U$. Given a function $u\in L^{1}(M,g)$, define the variation of $u$ by
      \begin{equation}
                 |Du|(M):=\sup\left\{\int_{M}u\,\mathrm{div}_g(X)\,\mathrm dv_g: X\in\mathfrak{X}_c(M), ||X||_{\infty}\leq 1\right\},
      \end{equation}  
      where $||X||_{\infty}:=\sup\left\{|X_p|_{g}: p\in M\right\}$ and $|X_p|_{g}$ is the norm of the vector $X_p$ in the metric $g$ on $T_pM$.
      We say that a function $u\in L^{1}(M,g)$, has \textbf{bounded variation}, if $|Du|(M)<\infty$ and we define the set of all functions of bounded variations on $M$ by $BV(M,g):=\big\{u\in L^1(M,g):\:|Du|(M)<+\infty\big\}$. A function $u \in L^1_{\text{loc}}(M)$ has \textbf{locally bounded variation in $M$}, if
for each open set $U \Subset M$, 
\[|Du|(U):=\sup\left\{ \int_U u\, \mathrm{div}_g(X)\,\mathrm dv_g : X\in\mathfrak X_c(U), \|X\|_{\infty}\le1  \right\}<\infty,\]
and we define the set of all functions of locally bounded variations on $M$ by $BV_{loc}(M) := \{u \in L^1_{loc}(M) : |Du|(U) < +\infty, U\Subset M\}$. So for any $u\in BV(M,g)$, we can associate a vector Radon measure on $M$, $\nabla^g u$ with total variation $|\nabla^gu|$.
      \end{definition}
\begin{definition}\label{Def:LocallyFinitePerimeterSets} Let $(M,g)$ be a Riemannian manifold of dimension $N$, $U\subseteq M$ be an open subset, $\mathfrak{X}_c(U)$ the set of smooth vector fields with compact support in $U$. Given $E\subset M$ measurable with respect to the Riemannian measure, the \textbf{perimeter of $E$ in $U$}, $ \mathcal{P}_g(E, U)\in [0,+\infty]$, is
      \begin{equation}
                 \mathcal{P}_g(E, U):=\sup\left\{\int_{U}\chi_E\mathrm{div}_g(X)dv_g: X\in\mathfrak{X}_c(U), ||X||_{\infty}\leq 1\right\},
      \end{equation}  
where $||X||_{\infty}:=\sup\left\{|X_p|_{g}: p\in M\right\}$ and $|X_p|_{g}$ is the norm of the vector $X_p$ in the metric $g$ on $T_pM$. If $\mathcal{P}_g(E, U)<+\infty$ for every open set $U\subset\subset M$, we call $E$ a \textbf{locally finite perimeter set}. Let us set $\mathcal{P}_g(E):=\mathcal{P}_g(E, M)$. Finally, if $\mathcal{P}_g(E)<+\infty$ we say that \textbf{$E$ is a set of finite perimeter}. We will use also the following notation $\mathcal{P}_g(E, F):=|\nabla\chi_E|_g(F)$ for every Borel set $F\subseteq M$. 
\end{definition}
\begin{proposition}[Riemannian version of \smash{\cite[Proposition~2, p.\ 133]{Modica}}]\label{prop:2Modica}
Let $(M^N,g)$ be a complete smooth Riemannian manifold, let $A$ and $\Omega$  be open subsets of $M$ with $\partial A$ a non-empty, compact, smooth hypersurface, and with $\mathcal{H}_g^ {N-1}(\partial A\cap\partial\Omega)=0$. Assume that $\overline\Omega$ is compact with smooth boundary (possibly empty). Given real numbers $\alpha,\beta$, with $\alpha<\beta$, define the function $v_0:\Omega\rightarrow\R$ by
\[
v_0(x) =\begin{cases} \alpha, & \text { if } x \in A,\\ \beta, &\text { if } x \in \Omega\setminus  A .\end{cases}
\]
Then there is a family $\left( v_\varepsilon\right)_{\varepsilon > 0 }$ of Lipschitz continuous function on $M$ such that $v_\varepsilon$ converges to $v_0$ in $L ^ { 1 } (\Omega)$ as $\varepsilon \rightarrow 0 ^ { + }, \alpha \leq v_\varepsilon\leq \beta$ for every $\varepsilon > 0 ,$ and 
\begin{enumerate}[(i)] 
\item $\int_\Omega v_\varepsilon\,\mathrm  dv_g= \int_\Omega v_0\,\mathrm dv_g=\alpha | A \cap \Omega | + \beta | \Omega \setminus  A |, \quad \forall \varepsilon > 0$,
\item\label{Eq:(ii)ModicaProp2}  $\limsup\limits_{\varepsilon\rightarrow0^+}E_{\varepsilon}(v_\varepsilon)\leq\mathcal{P}_g(A,\Omega)\,\sigma(\alpha,\beta)$, where $\sigma(\alpha,\beta)=\int_\alpha^\beta\sqrt{2W(s)}\,\mathrm ds$, and $E_\varepsilon$ is as in \eqref{Eq:DefinitionOfFunctional}.
\end{enumerate}
\end{proposition}

\begin{proof}
 Let us define the function $d_A$ as 
 $$
d_A ( x ) = \left\{ \begin{aligned} - \operatorname { dist } ( x , \partial A ) & \text { if } x \in A \\ \operatorname { dist } ( x , \partial A ) & \text { if } x \notin A. \end{aligned} \right.
$$
\begin{figure}[H]
\begin{center}
\includegraphics[]{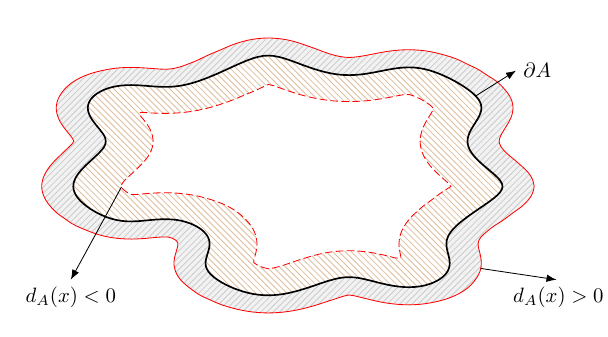}
\caption{Level sets of a typical example of the signed distance function for a subset $A\subseteq M$.}\end{center}
\end{figure}
It is well known (see for instance \cite[Lemma~4]{Modica} that $d_A$ is Lipschitz continuous, 
that $|\nabla_g d_A(x)|_g= 1$ for almost all $x\in M$ and that, if
$S_t:= \left\{ x \in M: d_A ( x ) = t \right\}$, then
\begin{equation}  
\lim _ { t \rightarrow 0 }\mathcal{H }_g^{N-1}\left( S _ { t } \cap \Omega \right) =\mathcal{H }_g^{N-1}( \partial A \cap \Omega ).
\end{equation}
Define the function $q _ { 0 } : \mathds  { R } \rightarrow \mathds  { R }$ by
\[
q _ { 0 } ( t ) = \alpha \quad \text { if } t < 0 , \quad q _ { 0 } ( t ) = \beta \quad \text { if } t \geq 0,
\]
and let $v_{0}( x ) = q_0({d_A(x)})$. 

Now, let us consider functions $q_\varepsilon$ satisfying the ordinary differential equation:
\begin{equation}\label{eq:ode1}
\frac{\varepsilon}2 q_ { \varepsilon } ^ { \prime  }(t) ^2= \frac1\varepsilon W \left( q_ { \varepsilon } (t)\right)+
\frac{\sqrt { \varepsilon }}2;
\end{equation}
these functions give an approximation of $q_0$, and will be employed to define the desired maps  $\left( v_ { \varepsilon } \right)$. 
Let us explain why this equation.
We want  approximate the two-valued function $q_{ 0 }$ by a Lipschitz continuous function $q_ { \varepsilon }$, which interpolates between $\alpha$ and $\beta$ and, at the same time, minimizes the 
one-dimensional Van der Walls-Allen-Cahn-Hilliard gradient phase field energy functional
\[
\int_{\R} \left[ \frac\varepsilon2 q_ { \varepsilon } ^ { \prime 2 } + \frac1\varepsilon W \left( q_ { \varepsilon } \right) \right] \mathrm d t.
\]
The corresponding Euler equation is
$\varepsilon^2 q_ { \varepsilon } ^ { \prime \prime } =  W ^ { \prime } \left( q_ { \varepsilon } \right)$; 
multiplying by $q_ { \varepsilon } ^ { \prime }$ and integrating, we obtain $\frac{\varepsilon^2}2q_\varepsilon ^ { \prime 2 } =  W \left( q_ { \varepsilon } \right)+c_ { \varepsilon } $.
To avoid the constant trivial solutions, the constant $c_\varepsilon$ cannot be set equal to zero.
On the other hand, we need $c _ { \varepsilon } \gg \varepsilon^2$ to make $q _ { \varepsilon }$ fill the gap between $\alpha$ and $\beta$ as quickly as possible (note that $q'^2_ { \varepsilon }\geq 2c _ { \varepsilon } / \varepsilon^2 )$), and for that reason we choose $c _ { \varepsilon } = \varepsilon^{3/2}/2$.
%
%
To construct the functions $q_\varepsilon$, consider for a fix $\varepsilon > 0$ the function
\[
\psi _ { \varepsilon } ( t ) = \int _ { \alpha } ^ { t }  \frac { \varepsilon}{ \sqrt{\varepsilon^{3/2}+2W(s)}}\,\mathrm ds, \quad \alpha \leq t \leq \beta
\]
where $\eta _ { \varepsilon } = \psi_ { \varepsilon } ( \beta )$.

Let $\widetilde q_ { \varepsilon } : \left[ 0 , \eta_ { \varepsilon } \right] \rightarrow [ \alpha , \beta ]$ denote the inverse of $\psi _ { \varepsilon }$, see Figure~\ref{fig:qepsilon}. 
\begin{figure}[H]
\begin{center}
\includegraphics[]{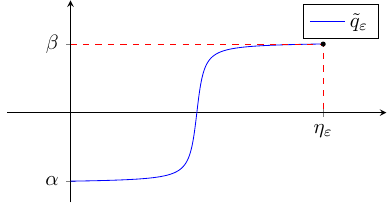}
\caption{Graph of $\widetilde q_\varepsilon$.}
\label{fig:qepsilon}
\end{center}
\end{figure}
Since $W$ is non-negative,
\begin{align}\label{eq:27}
0 < \eta_ { \varepsilon } \leq \varepsilon ^ { \frac { 1 } { 4 } } ( \beta - \alpha );
\end{align}
and, by the continuity of $W$, $\widetilde q_\varepsilon$ is of class $C ^ { 1 }$ and
\begin{align}\label{eq:28} 
\varepsilon \widetilde q_ { \varepsilon } ^ { \prime  }(t) = \sqrt{\varepsilon^{3/2}+2W(\widetilde q_\varepsilon)},
\end{align}
for $0 \leq t \leq \eta _ { \varepsilon }$.

We now extend the definition of $\widetilde q_ { \varepsilon }$ to the entire real line by setting 
$$
\widetilde q_ { \varepsilon } ( t ) = \alpha \quad \text { for } t < 0 , \quad \widetilde q_ { \varepsilon } ( t ) = \beta \quad \text { for } t > \eta_\varepsilon,
$$
so that $\widetilde q_ { \varepsilon }$ is a Lipschitz continuous function on $\R$. 
Note that, for every $t \in \mathds  { R }$, 
$\widetilde q_ { \varepsilon } ( t ) \leq q_ { 0 } ( t )$ and $\widetilde q_ { \varepsilon } \left( t + \eta_ { \varepsilon } \right) \geq q_ { 0 } ( t ) .$ Thus, there exists $\delta_{\varepsilon,A,V}\in\left[ 0 , \eta_ { \varepsilon } \right]$ such that
\begin{align}\label{eq:29}
\int_{\Omega } \widetilde  q_ { \varepsilon } \big( d_{A}(x)+\delta_{\varepsilon,A,V}\big) \mathrm dv_g = \int _ { \Omega }  q_ { 0 } (d_{A}( x ) ) dv_g= \int _ { \Omega } v _ { 0 } ( x ) \mathrm dv_g.
\end{align}
Finally, we define $q_ { \varepsilon } ( t ) = \widetilde q_ { \varepsilon } \left( t + \delta_{\varepsilon,A,V} \right)$ for $t \in \mathds  { R }$ and $$v_{\varepsilon,A,V}( x ) = q_ { \varepsilon } (d_{A}( x ) )=\widetilde q_ { \varepsilon } \left(d_{A}(x)+ \delta_ { \varepsilon,A,V} \right),$$ for
$x \in \Omega$.

We now prove that $v_{\varepsilon,A,V}\to v_0$ in $L^1$.
Notice that each $v_{\varepsilon,A,V}$ is a Lipschitz continuous function and $\alpha \leq v_{ \varepsilon,A,V } \leq \beta$.

By Lemma $4$ of \cite{Modica} and the co-area formula for Lipschitz functions (see for instance \cite{Federer})
\begin{align}\label{eq:30} 
\int _ { \Omega } f ( u ( x ) ) | \nabla_g u ( x ) |_g\,\mathrm dv_g= \int _ { \mathds  { R } } f ( t ) \mathcal{H}_g^{N-1}\big( \{ x \in \Omega : u ( x ) = t \}\big)\,\mathrm d t,
\end{align}
which holds for any Lebesgue measurable function $f$ and any Lipschitz continuous function $u$, we get the following
$$
\begin{aligned}
\int _ { \Omega } \left|v_{\varepsilon, A,V}-v_{0}\right|\mathrm dv_g
&= \int _ { \Omega } \left|  q_ { \varepsilon } ( d_A ( x ) ) -  q_ { 0 } ( d_A ( x ) ) \right| |\nabla_gd_A(x)|_g \mathrm dv_g\\
& = \int _ { -\delta_{\varepsilon,A,V}}^{\eta_\varepsilon-\delta_{\varepsilon,A,V}}\left| q_ { \varepsilon } ( t ) - q _ { 0 } ( t ) \right|\mathcal{H}_g^{N-1}\left( S _ { t } \cap \Omega \right)\mathrm  d t \\ 
& \leq\eta_{\varepsilon}(\beta - \alpha ) \sup _ { | t | \leq \eta _ { \varepsilon } } \mathcal { H }_g^{N-1} \left( S _ { t } \cap \Omega \right)\\
& \le \varepsilon ^ { \frac { 1 } { 4 }}  ( \beta - \alpha )^2 C(M,g, \Omega, A),\end{aligned}
$$
where $S _ { t } = \left\{ x \in M: d_A ( x ) = t \right\}$, and we obtain the last inequality applying \eqref{eq:27}. Then we conclude that $v_{\varepsilon,A,V}$,
converges to $v_{0}$ in $L ^ { 1 } ( \Omega )$ as $\varepsilon \rightarrow 0 ^ { + }$. 

To prove (ii) we call
$$
\gamma _ { \varepsilon,V} = \sup _ { | t | \leq \eta_ { \varepsilon } } \mathcal{H}_g^{N-1}\left( S _ { t } \cap \Omega \right).
$$
We again employ the coarea formula \eqref{eq:30}, obtaining
\[
\begin{aligned}
E_ {\varepsilon} \left(v_{\varepsilon, A,V}\right) 
& = \int _ { \mathds  { R } } \left[   \frac { \varepsilon } { 2 }q_{\varepsilon}^{\prime }(t)^2+ \frac 1{ \varepsilon }  W \left( q_ { \varepsilon } ( t ) \right) \right] \mathcal{H}_g^{N-1}\left( S _ { t } \cap \Omega \right) \mathrm d t\\ 
& \leq \gamma_ { \varepsilon } \int _ { -\delta_{\varepsilon,A,V}}^{\eta_\varepsilon-\delta_{\varepsilon,A,V}} \left[\frac { \varepsilon } { 2 } \widetilde q_{\varepsilon}^{\prime }(t+\delta_{\varepsilon,A,V})^2+ \frac 1{ \varepsilon }  W \left(\widetilde q_ { \varepsilon } ( t+\delta_{\varepsilon,A,V} ) \right) \right] \mathrm d t\\
&\leq \gamma _ { \varepsilon } \int _ { 0 } ^ { \eta_ { \varepsilon } } \left[\frac { \varepsilon } { 2 } \widetilde q_{\varepsilon}^{\prime }(t)^2+ \frac 1{ \varepsilon }  W \left(\widetilde q_ { \varepsilon } (t) \right) +\frac{\varepsilon^{1/2}}2\right] \mathrm d t
\end{aligned}
\]
and, recalling \eqref{eq:28},
$$
\begin{aligned}
E_ {\varepsilon} \left(v_{\varepsilon, A,V}\right) & \leq \gamma_ { \varepsilon } \int _ { 0 } ^ { \eta_ { \varepsilon } }  \left(  2W( \widetilde q_ { \varepsilon } ( t )) +\varepsilon^{3/2}  \right) ^ { \frac { 1 } { 2 } } \widetilde q _ { \varepsilon } ^ { \prime } ( t ) d t \\ 
& =\gamma_ { \varepsilon } \int _ { \alpha } ^ { \beta } ( 2W ( s )+\varepsilon^{3/2} ) ^ { \frac { 1 } { 2 } } d s. \end{aligned}
$$
Since Lemma $4$ of \cite{Modica} implies
$$
\lim _ { \varepsilon \rightarrow 0 ^ { + } } \gamma _ { \varepsilon } = \mathcal{H}_g^{N-1} ( \partial A \cap \Omega ) = \mathcal{P}_g( A ,\Omega),
$$
we conclude that
$$
\limsup _ { \varepsilon \rightarrow 0 ^ { + } }  E_ { \varepsilon } \left( v _ { \varepsilon,A,V } \right) \leq \mathcal{P}_g ( A ,\Omega) \int _ { \alpha } ^ { \beta } \sqrt{2W}d s=\mathcal{P}_g(A,\Omega)\sigma(\alpha,\beta),
$$
and the proposition is proved.
\end{proof}
\begin{definition}\label{thm:modicaapprox}
Given sets $A,\Omega\subset M$ and a function $v_0$ as in Proposition~\ref{prop:2Modica}, a family $(v_\varepsilon)_{\varepsilon>0}$ of functions as in the statement of the Proposition will be called a \emph{Modica approximation of $v_0$}. When the set $\Omega$ is not specified, it will be implicitly assumed $\Omega=M$.
\end{definition}

\section{The photography method}\label{sec:photography}
In this section we discuss a technique, originally due to Benci, Cerami, and others, (see \cite{BenciCeramiPassaseo} or \cite{BenciCeramiCalcVar}) 
which is a twist over the classical Lusternik-Schnirelmann theory and Morse theory. We will call this technique the \emph{photography method}, for reasons that will be clear along the way. 
A formal statement of the results generated by this method is given in Theorems \ref{Thm:FotoLS} and \ref{Thm:FotoMorse}.
\subsection{General setup} 
We start off by recalling a few basic definitions.
\begin{definition}\label{Def:Lusternik-SchnirelmannCategory} Let $(X,\tau)$ be a topological space and $Y\subseteq X$ be a closed subset. 
We define the \textbf{\emph{Lusternik-Schnirelmann category of $Y$ in $X$}} , denoted by $\cat_{X}(Y)$, as the minimum number $k\in\mathds {N}$ such that there exist open subsets $\mathcal{U}_{i}, \ldots, \mathcal{U}_{k} \subseteq X$ satisfying $Y\subseteq\bigcup_i\mathcal{U}_i$. If no such finite family exists, then one sets $\cat_X(Y)=+\infty$. Furthermore, 
one defines $\cat(X)=\cat_X(X)$.
\end{definition}
\begin{definition} Let $\mathfrak{M}$ be a $C^2$-Hilbert manifold, $J:\mathfrak{M}\to\R$ a $C^1$ functional, and $(u_n)$ a sequence in $\mathfrak{M}$. We say that $u_n$ is a \textbf{\emph{Palais--Smale sequence}} 
$($a \textbf{\emph{PS-sequence}} for short$)$, if 
\begin{equation}\label{Eq:DefPS}
J(u_n)\to c,
\end{equation}
\begin{equation}\label{Eq:DefPS0}
\big\Vert\mathrm d J(u_n)\big\Vert_{T_{u_n}^*\mathfrak{M}}\to 0.
\end{equation} 
\end{definition}
\begin{definition} Let $\mathfrak{M}$ be a $C^2$-Hilbert manifold, $J:\mathfrak{M}\to\R$ a $C^1$ functional. We say that \textbf{\emph{$J$ satisfies the Palais-Smale condition}}, if every Palais-Smale sequence has a convergent subsequence. 
\end{definition}
Classical results of Calculus of Variations relate the number of critical points in a sublevel of the energy functional with suitable topological/homological/cohomological invariants (category, Betti numbers, cuplength, etc.) of the sublevel. However, it is in principle rather involved to have a good topological description of subleveles of an abstract functional, which typically are the closure of arbitrary open subsets of infinite-dimensional manifolds.

The \emph{photography method} is a technique that allows us to estimate, when the functional space consists of real-valued functions on a manifold, the value of the topological invariants 
of the sublevels in terms of the analogous invariants associated to the underlying manifold. The estimate is obtained by \emph{reproducing} a copy of the underlying manifold in a given sublevel (the photography); this is done by means of a continuous function which associates to each point of the manifold, a map in the function spaces which concentrates its mass around the given point. The technique works when such operation can be made in such a way that the photography of the underlying manifold is sufficiently \emph{ample} in the sublevel, i.e., 
when the sublevel can be continuously retracted to the image of the photography. In many situations, such retraction is obtained as a \emph{barycenter map}.
This is formalized using continuous maps and homotopies, as follows.
\begin{theorem}\label{Thm:FotoLS} Let $\mathfrak{M}$ be a $C^2$-Hilbert manifold, $J:\mathfrak{M}\to\R$ be a $C^1$ functional, and  $J^c=\big\{v\in\mathfrak{M}:\:J(v)\le c\big\}$. Assume that 
\begin{enumerate}[(i)]
 \item\label{Assumption:Thm:FotoLS} $\inf\big\{J(u):u\in\mathfrak{M}\big\}>-\infty$,
 \item\label{Assumption:Thm:FotoLS0}  $J$ satisfies $(PS)$ condition, 
 \item\label{Assumption:Thm:FotoLS1}there exist $c>\inf J$, a topological space $X$ and two continuous maps $$f:X\to J^c,$$ $$g:J^c\to X$$ such that $g\circ f$ is homotopic to the identity map of $X$.
 \end{enumerate} 
 Then, there are at least $\cat(X)$ critical points in $J^c$. Furthermore, if $\mathfrak{M}$ is contractible and $\cat(X)>1$, there is at least one critical point $u\notin J^c$. 
\end{theorem}
\begin{proof} See \cite{BenciCeramiPassaseo} or \cite{BenciCeramiCalcVar}.
\end{proof}

The above result can be made slightly more accurate (at least in the nondegenerate case) by using Morse theory.
\begin{definition}\label{thm:defpoincarepolynomial}
Let $X$ be a topological space; the {\emph{Poincare's Polynomial $P_t(X)$ of $X$}} is defined as the following power series in the variable $t$ 
\begin{equation}\label{Eq:DefPP}
P_t(X)=\sum_{n=0}^\infty\beta_n(X)\,t^n.
\end{equation}
\end{definition}
\begin{remark} If $X$ is a compact manifold, we have that  $H^n(X)$ is a finite dimensional vector space and $H^n(X)$ is trivial for sufficiently large $n$. So the formal series \eqref{Eq:DefPP} is actually a polynomial.
\end{remark}


In the following definition, we give the notion of Morse index of a critical
point, which is necessary in our treatment to establish a relation between
the Poincare's polynomial $P_{t}(M)$ and the number of solutions to
the Euler equation associated to a given functional $J$. In this work we use the approach to Morse theory developed in \cite{BenciNewApproach}, which is suitable in problems arising from PDE's.  
\begin{definition}[Morse Index] Let $\mathfrak{M}$ be a $C^2$-Hilbert manifold, $J:\mathfrak{M}\to\R$ a $C^1$ functional
and let $u\in\mathfrak{M}$ an isolated critical point of $J$ at level\footnote{This means that $J(u)=c$, $\mathrm dJ(u)=0$, and there exists a neighborhood $\mathcal{U}$ of $u$ in $\mathfrak{M}$ such that the only critical critical point of $J$ contained in $\mathcal{U}$ coincide with $u$.} $c\in\R$. We denote by $i_t(u)$ the following formal power series in $t$ 
\begin{equation}
i_t(u)=\sum_{k=0}^{+\infty}\beta_k\big(J^c,J^c\setminus\{u\}\big)t^k.
\end{equation}
We call $i_t(u)$ the \textbf{\emph{(polynomial) Morse index of $u$}}. The number $i_1(u)$ is called the \textbf{\emph{multiplicity of $u$}}.  
\end{definition}%
If $J$ is of class $C^2$ in a neighborhood of $u$ and $J''[u]$ is not degenerate, we say that $u$ is nondegenerate. In this case we have that 
\begin{equation}\label{Eq:DefNondegenerateC^2}
i_t(u)=t^{\mu(u)},
\end{equation} where $\mu(u)$ is the \textbf{(numerical) Morse index of $u$}, i.e., the dimension of the maximal subspace on which the bilinear form $J''[u](\cdot, \cdot)$ is negative-definite. This fact suggests the following definition.
\begin{definition}\label{Eq:DefTopologicallyNondegenerate} Let $\mathfrak{M}$ be a $C^2$-Hilbert manifold, $J:\mathfrak{M}\to\R$ be a $C^1$ functional and let $u\in\mathfrak{M}$ be an isolated critical point of $J$ at level $c$. We say that $u$ is \textbf{\emph{$($topologically$)$ nondegenerate}}, if $i_t(u)=t^{\mu(u)}$, for some natural number $\mu(u)\in\mathds {N}$.
\end{definition}

\begin{theorem}\label{Thm:FotoMorse} Under assumptions \eqref{Assumption:Thm:FotoLS}, \eqref{Assumption:Thm:FotoLS0},  and \eqref{Assumption:Thm:FotoLS1} of Theorem \ref{Thm:FotoLS}, there exists $c_1>c$ such that one of the two following conditions hold: 
\begin{enumerate}
          \item $J^{c_1}$ contains infinitely many critical points.
          \item $J^{c}$ contains $P_1(X)$ critical points and $J^{c_1}\setminus J^{c}$, contains $P_1(X)-1$ critical points if counted with their multiplicity. More exactly we have the following relation
\begin{equation}\label{Eq:MorseRelation}
\sum_{u\in\Crit(J^{c_1})}i_t(u)=P_t(X)+t[P_t(X)-1]+(1+t)Q(t),
\end{equation}
where $Q(t)$ is a polynomial with nonnegative integer coefficients, and $\Crit(J^{c_1})$ denotes the set of critical points of $J$ in the sublevel $J^{c_1}$. 
\end{enumerate}
\end{theorem}
\noindent 
In particular, if all the critical points are nondegenerate there are at least $P_1(X)$ critical points with energy less or equal than $c$, and at least $P_1(X)\!-\!1$ with energy between $c$ and $c_1$.
\begin{remark} If we count the critical points with their multiplicity, then by Theorem \ref{Thm:FotoMorse} it follows that there are at least $2P_1(X)-1$ critical points. Namely, when the critical points are isolated, the result follows from the Morse's relation \eqref{Eq:MorseRelation},  otherwise there are infinitely many of them.
\end{remark}
\begin{remark}\label{thm:remsuperjacent}
Given topological spaces $X$ and $Y$, we say that \emph{$Y$ is homotopically superjacent to $X$} if there exist continuous maps
$f\colon X\to Y$ and $g\colon Y\to X$ such that $g\circ f$ is homotopic to the identity map of $X$. Thus, assumption \eqref{Assumption:Thm:FotoLS1} of Theorem~\ref{Thm:FotoLS} (and of Theorem~\ref{Thm:FotoMorse}) says that the sublevel $J^c$ is homotopically superjacent to $X$.
When $Y$ is homotopically superjacent to $X$, then the induced map $f_*$ in homotopy or homology is injective, which implies that $\cat(X)\le\cat(Y)$, and that for all $n\in\mathds N$,  $\beta_n(X)\le\beta_n(Y)$. This is the reason why the estimates on the number of critical points in 
Theorem~\ref{Thm:FotoLS} and in Theorem~\ref{Thm:FotoMorse} are given in terms of the topological invariants of $X$.
\end{remark}

\subsection{The photography method in our concrete setting}\label{sub:concrete}
In this section we will define the objects needed for the setup and the proof of Theorems~\ref{Thm:FotoLS} and \ref{Thm:FotoMorse}; an analysis of these objects will be carried out in the following sections. \smallskip

The objects $\mathfrak M$, $J$, $X$, $c$, $f\colon X\to J^c$ and $g\colon J^c\to X$ that appear in the statement of Theorem~\ref{Thm:FotoLS}
in our concrete setting are described below.
\begin{itemize}
\item $\mathfrak{M}=\mathfrak{M}^{V}$, where
\begin{equation*} 
\mathfrak{M}^{V}=\left\{ u\in H^{1}(M):\int_{M}u(x)\,\mathrm dv_g=V\right\}, 
\end{equation*}
\item $J={E_\varepsilon}\vert_{{\mathfrak{M}^V}}$, where 
\begin{equation*}
E_\varepsilon(u)=\frac\varepsilon2\int_{M} \left\vert
\nabla u\right\vert^2\,\mathrm dv_g+\frac{1}{\varepsilon}\int_{M} W(u(x))\,\mathrm dv_g,
\end{equation*}
\item $X=M$, $f=\Phi_{\varepsilon,V}:M\to E_{\varepsilon}^c\cap\mathfrak{M}^V=:\mathfrak{M}^V_{\varepsilon,c}$, where $c=c(\varepsilon,V)=\sigma c_NI_M(V)+\delta$ where 
\[\sigma=\sigma(0,1)=\int_0^1\sqrt{2W(s)}\,\mathrm ds,\] 
$\delta>0$ is a suitable small constant {that will be specified later 
(Lemma~\ref{Lemma:Concentration+})}, and 
$\Phi_{\varepsilon,V}\colon M\to \mathfrak M^V_{\varepsilon,c}$ is defined by:
\begin{equation}\label{eq:defphotmap}
\Phi_{\varepsilon,V}(x_0)(x):=U_{\varepsilon, V, x_0}\left(x\right),
\end{equation}
where $U_{\varepsilon,V,x_0}\colon M\to\mathds R$ is the function obtained in Proposition \ref{prop:2Modica} assuming $\alpha=0$, $\beta=1$, $\Omega:=M$, and $A:=M\setminus B_g(x_0, r_V)$ where $B_g(x_0, r_V)$ is the metric ball of volume $\mathrm{vol}_g(B_g(x_0, r_V))=V$. We observe that $U_{\varepsilon,V,x_0}\colon M\to\mathds R$ is a Lipschitz continuous function with the following properties:
\begin{itemize}
\item[$\bullet$] as it is easy to see, by construction we always have \[supp\left(U_{\varepsilon,V, x_0}\right)\Subset B_g(x_0,r_V+\delta_{\varepsilon, M\setminus B_g(x_0,r_V),V}),\] where $\delta_{\varepsilon, M\setminus B_g(x_0,r_V),V}>0$ is defined as in Proposition \ref{prop:2Modica}. Then 
for small $V\ll 1$ and small $\varepsilon\ll 1$ we have that for every $x_0$ it holds $$B_g(x_0,r_V+\delta_{\varepsilon, M\setminus B_g(x_0,r_V), V})\Subset B_g\left(x_0, \frac{\inj_M}2\right).$$
\item[$\bullet$] the family $\big(U_{\varepsilon,V,x_0}\big)_{\varepsilon>0}$ is a Modica approximation (see Definition~\ref{thm:modicaapprox}) of
the characteristic function of the geodesic ball $B_g(x_0,r_V)$ of volume $V$. Here, for $V\in\left]0,V_{x_0}\right[$ with:
\[C_1\le V_{x_0}=\vol_g\big(B_g(x_0,\tfrac12\inj_M)\big)\le {C_2(M,g)},\] 
\end{itemize}
where $C_1:=\vol_{g_{b^-}}(B_{(\mathbb{M}^N_{b^-},g_{b^-})}(0, \tfrac12\inj_M))$ and 
\[C_2:=\vol_{g_{b^+}}(B_{(\mathbb{M}^N_{b^+},g_{b^+})}(0, \tfrac12\inj_M)),\] $(\mathbb{M}^N_{k},g_{k})$ is the simply connected space form of constant sectional curvature $k\in\R$, and $b^-, b^+\in\R$ are such that $b^-\le Sec_g(\sigma, x)\le b^+$ for every $2$-dimensional subspace $\sigma\le T_xM$ and for every $x\in M$, with $Sec_g(\sigma, x)$ being the sectional curvature of the $2$-dimensional subspace $\sigma$ with respect to metric $g$. 
The existence of the functions $U_{x_0,V,\varepsilon}$ is proved in Proposition~\ref{prop:2Modica}.
\begin{figure}[H]
\begin{center}
\includegraphics[]{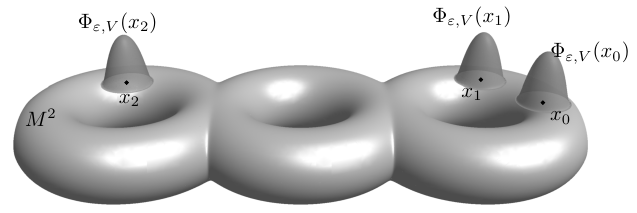}
\caption{The image of the photography map is a (photography) of the underlying (finite-dimensional) manifold in the infinite dimensional functional manifold $H^1(M)$.}\end{center}
\end{figure} 
\item assuming that the Riemannian manifold $M$ is isometrically embedded in some Euclidean space $\mathds R^l$ (Nash embedding theorem) 
$\hat{g}:=\pi\circ\beta:\mathfrak M^V_{\varepsilon,c}\to M$, where $\pi$ is the nearest point projection given in Definition~\ref{Def:TubularNeighborhood}, the barycenter map $\beta^* \colon \mathfrak M^V_{\varepsilon,c}\to\R^l$,
\begin{equation}\label{Eq:DefBarycenterNash1}
\beta^*(u):=\frac{\int_M xu(x)\,\mathrm dv_g(x)}{\int_M u(x)\,\mathrm dv_g(x)}.
\end{equation}
\end{itemize}
We will next show that  the above objects are well defined, and that they satisfy the assumptions required in the photography method.

\subsection{The Palais--Smale condition}\label{sub:PS}
First of all we need to prove the Palais--Smale condition for the functional ${E_{\varepsilon}}$ :
\begin{proposition}\label{Lemma:PS} For every $\varepsilon,V>0$, the functional ${E_{\varepsilon}}$ satisfies the Palais-Smale condition on $\mathfrak M^V$.
\end{proposition}

\begin{proof} Assume that $(u_n)$ is a Palais-Smale sequence for $E_\varepsilon$ in $\mathfrak M^V$; by density, we can suppose that $u_n$ is at least of class $C^2$ for all $n$. Observe that Equations \eqref{Eq:DefPS} and \eqref{Eq:DefPS0} written explicitly are
\begin{equation}\label{Eq:PSLS}
\frac{\varepsilon ^{2}}{2}\int_M \left\vert
\nabla u_n\right\vert ^{2}\,\mathrm dx+\int_M W(u_n(x))\,\mathrm dx\to c
\end{equation}
\begin{equation}\label{Eq:PSLS0}
-\varepsilon^2\Delta u_n+W'(u_n)=\lambda_n+T_n,
\end{equation}
where $(\lambda_n)_n$ is some sequence in $\mathds R$, and $T_n\to0$ strongly in $H^{-1}(M)$. 

Then, by \eqref{Eq:PSLS} and using the assumptions \eqref{Eq:Potential}, \eqref{Eq:Potential0} 
we have
\[
\begin{aligned} 
c+1 & \ge  \frac{\varepsilon ^{2}}{2}\int \left\vert\nabla u_n\right\vert ^{2}\,\mathrm dx+\int W(u_n(x))\,\mathrm dx\\
 & \ge   \frac{\varepsilon ^{2}}{2}\int\left\vert\nabla u_n\right\vert ^{2}\,\mathrm dx.
\end{aligned}
\]
Hence, $\left\vert\nabla u_n\right\vert$ is bounded in $L^2$. By Poincar\'e inequality on compact manifolds without boundary we know that there exists $C=C(M, g)>0$ such that $\left\|u_{n}-\bar{u}_{n}\right\|_{2}=\left\|u_{n}-V\right\|_{2} \leq C\left\|\nabla u_{n}\right\|_{2}$. 
Thus, $u_n$ is bounded in $H^1(M)$,
so that there exists  $u\in H^1(M)$ such that, up to a subsequence, $u_n\rightharpoonup u$, in $H^1(M)$.
Now the point is that we need to show that $u_n\to u$ strongly in $\mathfrak{M}_{\varepsilon,c}^{V}$.  

By \eqref{Eq:Potential0}, for some $p<\frac{2N}{N-2}$, the map $u\mapsto W'(u)$ of left composition with $W'$ gives a bounded nonlinear operator
from $L^p(M)$ to $L^q(M)$, with $\frac1p+\frac1q=1$; thus $q>\frac{2N}{N+2}\ge2$.
By the Sobolev embedding theorem, the inclusion $H^1(M)\hookrightarrow L^p(M)$ is compact, and thus
we get a compact nonlinear operator $H^1(M)\ni u\mapsto W'(u)\in L^q(M)$. This implies that, up to subsequences 
$W'(u_n)\to W'(u)$ strongly in $L^q(M)\subset H^{-1}(M)$. Multiplying \eqref{Eq:PSLS0} by $u_n$, integrating by parts the corresponding identity, and using the constraint $\int u_n=V$, we get that $\lambda_n$ is a bounded sequence. Whence, up to a subsequence we can assume $\lambda_n\to\lambda$. $\Delta:H^{1}(M)\to{H^{-1}(M)}$ is an isomorphism onto its image when restricted  to the subspace of functions orthogonal to the constants. Note that taking $c=-\frac V{\mathrm{vol}_g(M)}$ we get that $u_n+c$ is orthogonal to the space of constant functions and $\Delta(u_n+c)=\Delta u_n$ is $H^{-1}(M)$ convergent and therefore $u_n+c$ is strongly convergent in $H^1(M)$ and so $u_n$ is also strongly convergent $H^1(M)$.  
\end{proof}

\subsection{Continuity of the photography map}
This is the map $f$ that reproduces a copy of the finite-dimensional compact ambient manifold $M^N$ inside the infinite functional space which is the domain of the energy functional.
For the definition of $f$, see Section~\ref{sub:concrete}, formula \eqref{eq:defphotmap}.

Let us start by looking more closely at its definition and by proving the continuity of the photography map stated in the following proposition which to be stated needs the classical definition just below.
\begin{definition}\label{Def:IsPStrong}
The \textbf{\emph{isoperimetric profile function of $(M^N,g)$}} $($or briefly, the isoperimetric profile$)$ $I_{(M,g)}:\left[0,V(M)\right[\rightarrow\left [0,+\infty\right [$, is defined by $$I_M(V):= \inf\left\{A_g(\partial \Omega): \Omega\in \tau_M, \mathrm{vol}_g(\Omega )=V\right\},$$ where $\tau_M$ denotes the set of relatively compact open subsets of $M$ with smooth boundary, where $A_g$ is the $(N-1)$-volume form of $\partial\Omega$ induced by $g$.
\end{definition}
\begin{proposition}\label{Lemma:EnergyUpperBounds} For every $V\in\left]0, \mathrm{vol}_g(M)\right[$ and for every $\delta>0$ there exists $\varepsilon_1(V,\delta)>0$, such that for every $\varepsilon\in\left]0,\varepsilon_1\right[$ we have that $\Phi_{\varepsilon,V}$ carries $M$ into the sublevel $\mathfrak M^V_{\varepsilon,c}$,  where 
$c=\sigma I_M(V)+\delta$, and $\Phi_{\varepsilon,V}\colon M\to \mathfrak M^V_{\varepsilon,c}$ is a continuous function.  
\end{proposition}
\begin{proof} Recall from Section~\ref{sub:concrete} the map $\Phi_{\varepsilon,V}$ at some point $x_0\in M$ is defined in terms of Modica approximations for the characteristic functions of balls centered at $x_0$ with volume equal to $V$, see formula \eqref{eq:defphotmap}.
By \eqref{Eq:(ii)ModicaProp2} in Proposition~\ref{prop:2Modica} and the asymptotic expansion for small volumes of the area of the geodesic balls with respect to the enclosed volume, it follows that $E_{\varepsilon}(\Phi_{\varepsilon,V}(x_0))\lesssim\sigma I_M(V)$ as $\varepsilon\to0$, uniformly with respect to $x_0$ and $V$, where $I_M(V)$ is defined below in Definition \ref{Def:IsPStrong}. Thus using this and the compactness of $M$, 
the proposition is proved if we show that $\Phi_{\varepsilon,V}$ is continuous. 
To this aim, we will first prove the following estimate:
\begin{multline}\label{Eq:ContinuityEstimate}
||\Phi_{\varepsilon,V}(x_0)-\Phi_{\varepsilon,V}(x_1)||_{W^{1,2}(M)}  \le  C\big[||h_{x_0}-h_{x_1}||_{\infty}+|\delta_{\varepsilon,x_0,V}-\delta_{\varepsilon,x_1,V}|\big]\\
+ C\big[||\nabla h_{x_0}-\nabla h_{x_1}||_{\infty}\big],
\end{multline}
where $h_{x}(\cdot)=d_{\partial B_{g}\left(x, r_{V}(x)\right)}(\cdot)$ (see Figure~\ref{fig:hx}), $C=C(\varepsilon,V,M,g, W\big\vert_{[0,1]})>0$ and $\delta_{\varepsilon,x_0,V}:=\delta_{\varepsilon, M\setminus B_g(x_0,r_V), V}$, where $B_g(x_0,r_V)$ is the small geodesic ball enclosing volume $V$, see formula~\eqref{eq:29}. It is worth to notice here that for small volumes $V\ll 1$, we have that  $\partial B_g(x_0,r_V)$ is smooth. The desired continuity property of $\Phi_{\varepsilon,V}$ will follow from this inequality,
\begin{figure}[H]
\begin{center}
\includegraphics[]{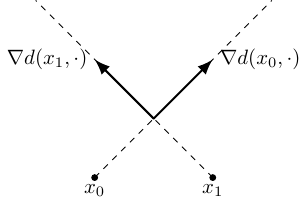}
\caption{Geometric illustration of the continuity of the photography map, because $\nabla d_{\partial B_{g}\left(x, r_{V}(x)\right)}(\cdot)=-\nabla d_{g}(x, \cdot)$.}
\label{fig:hx}
\end{center}
\end{figure} 
observing that:
\begin{itemize}
\item $||h_{x_0}-h_{x_1}||_{W^{1,\infty}(M)}\to0$, as $x_1\to x_0$;
\item $x\mapsto\delta_{\varepsilon,x,V}$ is a $C^1$ map, as it can be seen easily applying the implicit function theorem in \eqref{eq:29}.
\end{itemize}
Thus, $|\delta_{\varepsilon,x_0,V}-\delta_{\varepsilon,x_1,V}|\to0$ as $x_0\to x_1$ for any fixed $\varepsilon,V>0$.
In order to prove \eqref{Eq:ContinuityEstimate}, we proceed as follows:
\[
\begin{aligned}
||\Phi_{\varepsilon,V}(x_0)&-\Phi_{\varepsilon,V}(x_1)||_{W^{1,2}(M)}^2 
=\\ \int_M&|\tilde q'_\varepsilon(h_{x_0}(x)+\delta_{\varepsilon,x_0,V})\nabla h_{x_0}(x)-\tilde q'_\varepsilon(h_{x_1}(x)+\delta_{\varepsilon,x_1,V})\nabla h_{x_1}(x)|^2\,\mathrm dv_g\\
&\quad +  \int_M|\tilde q_\varepsilon(h_{x_0}(x)+\delta_{\varepsilon,x_0,V})-\tilde q_\varepsilon(h_{x_1}(x)+\delta_{\varepsilon,x_1,V})|^2\,\mathrm dv_g\\
& \le  2\|\tilde q_\varepsilon'\|_\infty^2\int|\nabla h_{x_0}(x)-\nabla h_{x_1}(x)|^2\,\mathrm dv_g\\
&+ 2\|\tilde q_\varepsilon'\|_\infty^2\int (| h_{x_0}(x)- h_{x_1}(x)|^2 +|\delta _ { \varepsilon ,  x _ { 1 } , V }-\delta _ { \varepsilon ,  x _ { 0 } , V  }|^2)\,\mathrm dv_g\\
&\stackrel{\eqref{eq:28}}{\le}C\left( \|\nabla h_{x_0}-\nabla h_{x_1}\|^2_\infty +\|h_{x_0}-h_{x_1}\|^2_\infty+|\delta _ { \varepsilon ,  x _ { 1 } , V  }-\delta _ { \varepsilon,  x _ { 0 } , V  }|^2 \right),
\end{aligned}
\]

where 
\begin{eqnarray*}
C=C(\varepsilon,V,M,g, W\vert_{[0,1]}) & = & 2\|\tilde q_\varepsilon'\|_\infty^2\mathrm{vol}_g(B_g(\cdot, \frac{\inj_M}2))\\
& \ge & C^{**}(\varepsilon, W\vert_{[0,1]})C^*(M,g)>0, 
\end{eqnarray*}
since $2\|\tilde q_\varepsilon'\|_\infty^2=C^{**}(\varepsilon, W\vert_{[0,1]})>0$, being $\tilde q'_\varepsilon$ the solution of the one-dimensional problem. From this estimate the continuity of $\Phi_{\varepsilon,V}$ follows easily and the theorem is proved.
\end{proof}
\subsection{The barycenter map}\label{sub:barycenter}
In this section we will show that the baricenter map \eqref{Eq:DefBarycenterNash1} is well-defined and continuous. We start with the following:
\begin{definition}\label{Def:TubularNeighborhood} Given an isometric embedding $i:(M^N,g)\to(\R^l,\xi)$. We define the \textbf{\emph{normal injectivity radius $r_{i}(M)$}} as the largest nonnegative number $r$ such that the normal exponential $exp_{\nu M}:\nu M\to\R^l$ is a diffeomorphism of a neighborhood of the zero section of $\nu M$ into $M_{r}$ where $M_r:=\{x\in\R^l: d_\xi(x,M)<r\}\subseteq\R^l$ and  $\nu M$ denotes the normal bundle induced by $i$ on $M$. Let us denote by $\pi:M_{r_i(M)}\to M$ the canonical projection associated with the canonical projection $\widetilde\pi:\nu M\to M$. 
\end{definition}
\begin{remark} Notice that $M$ is a retract of $M_{r_i(M)}$, and $r_i(M)>0$, since $M$ is compact. 
\end{remark}
For the reader's convenience, we give a proof of the following simple result.
\begin{lemma}\label{Lemma:BarycenterContinuity} The map $\beta^*:H^1(M)\setminus\{0\}\to\R^l$ defined in \eqref{Eq:DefBarycenterNash1} is continuous. 
In particular, their restrictions to $\mathfrak{M}^V$  are continuous for every $V\in\R$.
\end{lemma}
\begin{proof} Let us prove the continuity of $\beta^*$.

 For all $w\in H^1(M)$, set $\mu_w:=\int_M w(x)\,\mathrm dv_g(x)$. We have the following estimate
\begin{eqnarray}\label{Eq:BarycenterContinuity}
\left|\frac{\int_M xu(x)\,\mathrm dv_g(x)}{\int_M u(x)\,\mathrm dv_g(x)}-\frac{\int_M xv(x)\,\mathrm dv_g(x)}{\int_M v(x)\,\mathrm dv_g(x)}\right| & \le & \frac{||x||_\infty}{\mu_u}\int_M\left |u-\frac{\mu_u}{\mu_v}v\right|\,\mathrm dv_g,
\end{eqnarray}
where $||x||_{\infty}:=\sup_{x\in M}\{|x|_{\R^l}\}=C(i)<+\infty$, because $M$ is compact. Here $|x|_{\R^l}$ is the Euclidean length of the position vector and $i$ is the isometric embedding of $M$ in $\R^l$. It is easy to show that the right-hand side of 
\eqref{Eq:BarycenterContinuity} goes to zero when $v\to u$ in  $L^2(M)$ (Lebesgue's dominated convergence, H\"older inequality). 
\end{proof}

\begin{remark}\label{thm:remrangebarycenter}
In order to apply the abstract theory of Theorems \ref{Thm:FotoLS}, \ref{Thm:FotoMorse}, in our concrete setting, a crucial point to be shown is that 
for fixed small $\varepsilon,V>0$ and for $c$ close to the minimum of $E_\varepsilon$ in $\mathfrak M^V$, the image $\beta^*(\mathfrak M^V_{\varepsilon,c})$ is contained in a tubular neighborhood $M_r$ of $M$ in $\mathds R^l$, whose thinckness $r>0$ is small enough to make the nearest point projection $M_r\to M$ well-defined and continuous.
The proof of this fact is rather involved, and it requires notions and nontrivial results about the isoperimetric problem in Riemannian manifolds.
\end{remark}

We recall here a very classical notion of measure theory that will be useful in the proof of Proposition \ref{Prop:1and3Modica}.
\begin{definition} Let $u,(u_n)_n$ be measurable functions on a measure space $(X,\Sigma,\mu)$. The sequence $u_n$ is said to \textbf{\emph{converge globally in measure to $u$}}, if for every $\varepsilon>0$, it holds
\[\lim _{{n\to \infty }}\mu\big\{x\in X:|u(x)-u_{n}(x)|\geq \varepsilon \big\}=0.\]
\end{definition}
We need a Riemannian formulation of another result from \cite{Modica}. It is worth pointing out that it is natural that Modica's compactness argument can be carried out in the Riemannian setting. We also mention that Riemannian versions of the proof have appeared in the literature. For instance, in the aforementioned result of A. Pisante and F. Punzo \cite{Pisante2016}. Compare also with the one in A. Pisante and M. Ponsiglione \cite{Pisante2010}, where the proof is carried out for specific choices of nonlinearities $W$, in the hyperbolic space, whereas following the same strategy. 
\begin{proposition}[Riemannian version of \smash{\cite[Prop.~1, Prop.~3]{Modica}}]\label{Prop:1and3Modica} 
Under assumption (d) in Section~\ref{sec:formulation} for the potential $W$ (see \eqref{Eq:aPrioriL1Convergence+}),
assume also that there exist constants $E^*>0$, $t_0>0$, $0<c_1<c_2$, $2<p_1<\hat2$, $p_1\le p_2\le 2(p_1-1)$, with $\hat2:=\frac{2^*}2+1$, a sequence of positive numbers such that $\varepsilon_i\to0^+$, and a sequence of functions $u_{\varepsilon_i}\in H^1(M)$ satisfying 
\begin{equation}\label{Eq:AssumptionB3+}
E_{\varepsilon_i}(u_{\varepsilon_i})\le E^*,\forall i\in\mathds {N}.
\end{equation}
Then, there exists a subsequence still denoted $(\varepsilon_i)_i$  such that $(u_{\varepsilon_i})_i$ converges to a function $u_\infty\in BV(M)$ in $L^1(M)$. Moreover, there exists a finite perimeter set $\Omega$ such that $u_\infty=\chi_\Omega$ such that $\mathcal{P}_g(\Omega)=|Du_{\infty}|(M)=\int_M|Du_{\infty}|\le\frac{E^*}{\sigma}$, where $\sigma=\int_0^1\sqrt{2W(s)}ds$. In particular, $\int_{M} u_{\varepsilon_{i}} d v_{g} \rightarrow \operatorname{vol}_{g}(\Omega)$.
\end{proposition}
\begin{proof} Let $\phi$ be the primitive function of $\left(2W\right)^{\frac12}$ with $\phi(0)=0$, i.e., \[\phi(t)=\int_0^t\left(2W(s)\right)^{\frac12}\,\mathrm ds,\] and set $v_ {\varepsilon_i}(x):=\phi\left( u_ {\varepsilon_i}(x)\right)$. We claim that the family $\left(v_{\varepsilon_i}\right)_{\varepsilon_i>0}$ is bounded in $L^{1}(M)$. In fact if \eqref{Eq:aPrioriL1Convergence+} holds, it is not restrictive to assume that $t _ { 0 } \geq 1$, and we easily have that
\begin{eqnarray*}
\phi(t) & = & \int_0^ {t_0}\left(2W(s)\right)^{\frac12}ds + \int_{t_0}^t\left(2W(s)\right)^{\frac12}\,\mathrm ds\\
& \le & \int_0^ {t_0}\left(2W(s)\right)^{\frac12}ds+2\frac {\sqrt {2c _ { 2 } } } {p_2+2} t ^ {\frac{p_2}2+ 1 },\quad \forall t \geq t _ { 0 }.
\end{eqnarray*}
Moreover, $p_2\le2(p_1-1)$ implies that $\frac{p_2}2+1\leq p_1;$ hence
$$
\phi ( t ) \leq c' _ { 3 } + c' _ { 4 } W ( t ), \quad \forall t \geq 0,
$$
for some real constants $c'_3$ and $c'_4$. One can prove an analogous estimate for $t\le 0$, so we get 
$$
|\phi(t)|\leq c _ { 3 } + c _ { 4 } W ( t ), \quad \forall t\in\R,
$$
for some real constants $c_3$ and $c_4$.
Then,
$$
\begin{aligned} \int_M |v_ {\varepsilon_i}| dv_g& \leq c _ { 3 }\mathrm{vol}_g(M)+ c _ { 4 } \int_M W \left( u_{\varepsilon_i} ( x ) \right) \,\mathrm dv_g\\ & \leq c _ { 3 }\mathrm{vol}_g(M) + c _ { 4 }\varepsilon_iE_ {\varepsilon_i} \left( u_{\varepsilon_i} \right), \end{aligned}
$$
thus
\begin{eqnarray*}
\int_M |v_{\varepsilon_i}|dv_g\le c_3\mathrm{vol}_g(M) + \widetilde c_4E^*, \forall i \in\mathds {N},
\end{eqnarray*}
for some real constant $\widetilde c_4$, and from this we conclude that $ \left(v_{\varepsilon_i}\right)$ is a bounded sequence in $L^1(M)$. Since the functions $u_{\varepsilon_{i}}$ (and hence $v_{\varepsilon_{i}}$ ) have better regularity, by the chain rule we obtain easily that 
$$|\nabla v_{\varepsilon_i}|=|\phi'(u_{\varepsilon_i})\nabla u_{\varepsilon_i}|=\left(2W(u_{\varepsilon_i})\right)^{\frac12}|\nabla u_{\varepsilon_i}|.$$
From the elementary inequality $ab\le\frac{\eta a^2}2+\frac{b^2}{2\eta}$ valid for every $\eta>0$ and $a,b\in\R$ but nontrivial only when $a\cdot b\in\left]0,+\infty\right[$, putting $\eta:=\varepsilon_i$, $a=|\nabla u_{\varepsilon_i}|$, $b=\sqrt{2W(u_{\varepsilon_i})}$ we get 
\begin{eqnarray}\label{Eq:Modica13}\nonumber
\int_M|\nabla v_{\varepsilon_i}|\,\mathrm dv_g & \le & \int_M\left(\frac12\varepsilon_i|\nabla u_{\varepsilon_i}|^2+\frac1{\varepsilon_i} W(u_{\varepsilon_i})\right)\,\mathrm dv_g\\
& \le & E_{\varepsilon_i}(u_{\varepsilon_i})\stackrel{\eqref{Eq:AssumptionB3+}}{\le}E^*.
\end{eqnarray}
Applying the compactness theorem for bounded variation functions, (cf.\ \cite[Theorem~3.23]{AmbrosioFuscoPallara} or \cite[Proposition~1.4]{MPPP2007}), there exists a subsequence also denoted by $(v_{\varepsilon_i})_i$ and an a.e.\ pointwise limit function $v_{\infty}\in BV(M,g)$ that is the $L^1(M,g)$ limit of the $v_{\varepsilon_i}$, which satisfies  
\[|Dv_\infty|(M)\le\liminf_{i\to\infty}\left\|\nabla v_{\varepsilon_{i}}\right\|_{L^{1}(M)}\le E^*,\]
where as customarily we denote by $Dv_\infty$ the Radon measure representing the distributional derivative of $v_\infty$ and by $|Dv_\infty|$ is its total variation. We now return to the study of the original functions $u_{\varepsilon_i}$. Since our potential $W\in C^2(\R)$ is nonnegative and vanishes only at $0$ and $1$, we deduce that $\phi$ is at least $C^1$ (actually $\phi$ is $C^3$) and a strictly monotone increasing function. Let $\psi$ be the inverse function of $\phi$ which always exists and it is strictly monotone increasing as well. By the implicit function theorem we get that $\psi$ inherits the same regularity of $\phi$. Let us define $u_\infty( x ) = \psi \left(v_\infty( x )\right)$, from the chain rule for $BV$ functions we also obtain that $u_\infty\in BV(M)$.
By \eqref{Eq:aPrioriL1Convergence+} then $\phi ^ { \prime } ( t )\geq\sqrt {2c_1} t _ { 0 } ^ { p_1 / 2 }$ for every $|t|\geq t _ { 0 } ;$ hence $\psi$ is Lipschitz continuous on $]-\infty, \phi\left(-t_0\right)]\mathring{\bigcup}[ \phi \left( t _ { 0 } \right) , + \infty [$ and so uniformly continuous on the entire real line. From this combined with Theorem $2$ of \cite{BartleJoichi} we infer that $u _ {\varepsilon_i} = \psi \circ v_{\varepsilon _i} $ converges in measure on $M$ to $u_\infty$ as $\varepsilon_i\rightarrow 0^+$ so \emph{a fortiori} also  $u _ {\varepsilon_i}$ converges pointwise a.e. on $M$ to $u_\infty$ as $\varepsilon_i\rightarrow 0^+$; since
$$
\begin{aligned} 
\int_M |u _ {\varepsilon_i}|^{p_1}dv_g & \leq \int_M t_0^{p_1}dv_g+ \frac { 1 } { c _ { 1 } } \int _M W \left( u _ {\varepsilon _i} ( x ) \right) dv_g\\ 
&\leq  t_0^{p_1}\Vol_g(M) + \frac1{c_1} \varepsilon_iE_{\varepsilon_i} \left( u _ {\varepsilon_i} \right)\\
&\le t_0^{p_1}\Vol_g(M)+\frac{\varepsilon_i}{c_1}E^*,
\end{aligned}
$$
we conclude that $\left(u_{\varepsilon_i}\right)$ is bounded in $L^{p_1}(M)$ with $p_1\geq 2$. This implies (via H\"older inequality) uniform integrability of the sequence $\left(u_{\varepsilon_i}\right)_i$. Hence, by the classical theorem of Vitali (compare Theorem $2.18$ of \cite{AmbrosioDaPratoMennucci} or Theorem $4.5.4$ of  \cite{Bogachev2007measure}), we know that uniform integrability and convergence in measure (which implies pointwise convergence a.e.) that $\left(u_{\varepsilon_i} \right)$ actually converges in $L ^1(M)$ to $u_\infty$. The remaining part of the proof goes along the same lines of the proof of \cite[Proposition~1]{Modica}. In fact, by Fatou's Lemma and \eqref{Eq:AssumptionB3+} it holds 
\begin{eqnarray*}
0\le\int_MW(u_\infty)dv_g & \le & \liminf_{i\to\infty}\int_MW(u_{\varepsilon_i})dv_g\le\liminf_{i\to\infty}\varepsilon_iE_{\varepsilon_i}\\ & \stackrel{\eqref{Eq:AssumptionB3+}}{\le} & \liminf_{i\to\infty}\varepsilon_iE^*=0.
\end{eqnarray*}
The last chain of inequalities shows that $W(u_\infty)=0$ a.e. on $M$ which in turn implies that $u_\infty(M)=\{0,1\}$ a.e. As already observed above $u_\infty\in BV(M)$ and so the sets $u_\infty^{-1}(0)$ and $u_\infty^{-1}(1)$ are of finite perimeter in $M$, moreover we can apply Theorem $3.96$ of \cite{AmbrosioFuscoPallara} and formula $(3.90)$ of \cite{AmbrosioFuscoPallara} to get by the Fleming-Rishel formula
\begin{equation}
\mathcal{H}_g^{N-1}(\partial^*\{u_\infty^{-1}(1)\})=|Du_\infty|(M)=\int_M|Du_\infty|=\frac1{\sigma}\int_M|Dv_\infty|\stackrel{\eqref{Eq:Modica13}}{\le}\frac{E^*}{\sigma},
\end{equation}
where $\sigma=\sigma(0,1)$.
 This yields the proof of Proposition \ref{Prop:1and3Modica}.
\end{proof}
\begin{lemma}\label{Lemma:Concentration-} Let $(M^N,g)$ be a compact Riemannian manifold. For every $1>\eta>0$, $V\in\left]0, \mathrm{vol}_g(M)\right[$, $\delta>0$ there exists $\varepsilon_0=\varepsilon_0(g, W, \eta, V, \delta)>0$ such that for every $0<\varepsilon<\varepsilon_0$ and for any $u\in \mathfrak{M}^V_{\varepsilon,c}$ with $c=c(W, V,\delta)=\sigma I_M(V)+\delta>0$, there exists $\Omega_{V,u}$ a finite perimeter set of volume $V$ such that $||u-\chi_{\Omega_{V,u}}||_{L^1(M)}\le\eta$ and $\mathcal{P}_g(\Omega_{V,u})\le\frac{c}{\sigma}$.  
\end{lemma}
\begin{proof} We argue by contradiction. Suppose that the conclusion does not hold. Then there exist $1>\eta>0$, $V\in\left]0, \mathrm{vol}_g(M)\right[$, $\delta>0$ a sequence $\varepsilon_i\to 0$, $u_{\varepsilon_i}\in \mathfrak M^V_{\varepsilon_i,c}$ such that for every $\Omega_V$ finite perimeter set of volume $V$ we have 
\begin{equation}\label{Eq:ContradictionConcentration0}
||u_{\varepsilon_i}-\chi_{\Omega_V}||_{L^1(M)}>\eta>0.
\end{equation}
 If we assume furthermore that holds \eqref{Eq:aPrioriL1Convergence+} we can apply Proposition \ref{Prop:1and3Modica} with $E^*:=c$. This provides a subsequence still denoted $(\varepsilon_i)_i$, a finite perimeter set $\Omega_V$ of volume $V$ such that $\mathcal{P}_g(\Omega_{V,(u_{\varepsilon_i})_i})\le\frac{c}{\sigma}$ and $$||u_{\varepsilon_i}-\chi_{\Omega_{V,(u_{\varepsilon_i})_i}}||_{L^1(M)}\to0,\text{ as } i\to+\infty.$$   This last equation contradicts \eqref{Eq:ContradictionConcentration0} and in turn completes the proof of the lemma.
\end{proof}
Observe that we can choose $\delta$ sufficiently small and refine the result of Lemma~\ref{Lemma:Concentration-} in order to 
have that $\Omega_{u,V}$ as above is actually an isoperimetric region; this yields the following concentration lemma for functions with energy close to the minimum energy level.
\begin{lemma}\label{Lemma:ConcentrationIsoperimetric} Let $(M^N,g)$ be a compact Riemannian manifold. For every $1>\eta>0$, $V\in\left]0, \mathrm{vol}_g(M)\right[$, there exist $\delta_0=\delta_0(\eta,V, M^N, g, W)>0$ such that for every $0<\delta<\delta_0$ there exists $\varepsilon_0=\varepsilon_0(g, W, \eta, V, \delta)>0$ such that for every $0<\varepsilon<\varepsilon_0$ and for any $u\in \mathfrak M^V_{\varepsilon,c}$ with $c=c(W, V,\delta)=\sigma I_M(V)+\delta$ there exists $\Omega_{V,u}$ isoperimetric region of volume $V$ such that $||u-\chi_{\Omega_{V,u}}||_{L^1(M)}\le\eta$.  
\end{lemma}
\begin{proof} We argue by contradiction. Suppose that the conclusion does not hold. Then there exist $1>\eta>0$, $V\in\left]0, \mathrm{vol}_g(M)\right[$, a sequence $\delta_i\to0^+$ a second sequence $\varepsilon_i\to 0^+$, a sequence of functions $(u_{\varepsilon_i})_i$ satisfying $u_{\varepsilon_i}\in E_{\varepsilon_i}^{c_i=c(V,\delta_i)=\sigma I_M(V)+\delta_i}\cap\mathfrak M^V$ such that for every $\Omega_V$ isoperimetric region of volume $V$ we have 
\begin{equation}\label{Eq:ContradictionConcentration}
||u_{\varepsilon_i}-\chi_{\Omega_V}||_{L^1(M)}>\eta>0.
\end{equation}
If we assume furthermore that \eqref{Eq:aPrioriL1Convergence+} holds, we can apply Proposition~\ref{Prop:1and3Modica} with $E^*:=c_1$. 
This provides a subsequence, still denoted by $(\varepsilon_i)_i$, a finite perimeter set $\Omega^{(1)}_V$ of volume $V$ such that $\mathcal{P}_g(\Omega_V^{(1)})\le\frac{c_1}{\sigma}$ and 
\begin{equation}\label{Eq:ContradictionConcentration1+}
||u_{\varepsilon_i}-\chi_{\Omega^{(1)}_V}||_{L^1(M)}\longrightarrow0,\text{ as } i\to+\infty.
\end{equation}  
To this subsequence we apply again Proposition \ref{Prop:1and3Modica}, now with $E^*:=c_2$. In this way we obtain again a new subsequence, still denoted $(\varepsilon_i)_i$, a finite perimeter set $\Omega^{(2)}_V$ of volume $V$ such that $\mathcal{P}_g(\Omega_V^{(2)})\le{c_2}/{\sigma}$ and 
\begin{equation}\label{Eq:ContradictionConcentration2+}
||u_{\varepsilon_i}-\chi_{\Omega^{(2)}_V}||_{L^1(M)}\to0,\text{ as } i\to+\infty.
\end{equation} 
The sequence appearing in \eqref{Eq:ContradictionConcentration2+} being a subsequence of the sequence appearing in \eqref{Eq:ContradictionConcentration1+} readily gives that $\Omega_V^{(2)}=\Omega_V^{(1)}$ by the uniqueness of the limit. Continuing this process and applying a standard diagonal argument we get the existence of a subsequence still denoted $(\varepsilon_i)_i$, a finite perimeter set $\Omega_V=\Omega_V^{(1)}=\Omega_V^{(2)}=\Omega_V^{(3)}=...$, of volume $V$ such that 
\begin{equation}\label{Eq:ContradictionConcentration1-}
\mathcal{P}_g(\Omega_V)\le\frac{c_i}{\sigma},\forall i\in\mathds {N},
\end{equation}
and 
\begin{equation}\label{Eq:ContradictionConcentration3+}
||u_{\varepsilon_i}-\chi_{\Omega_V}||_{L^1(M)}\to0,\text{ as } i\to+\infty.
\end{equation}  
From \eqref{Eq:ContradictionConcentration1-} we conclude immediately that $\mathcal{P}_g(\Omega_V)\le I_M(V)$ and so a fortiori we can assert that $\Omega_V$ is an isoperimetric region of volume $V$. This last fact combined with equation \eqref{Eq:ContradictionConcentration3+} contradicts \eqref{Eq:ContradictionConcentration} and in turn completes the proof of the lemma.  
\end{proof}
\begin{lemma}[\smash{\cite[Theorem~2.2]{MorganJohnson2000}} and \smash{\cite[Lemma~4.9]{NardulliOsorioIMRN}}]
\label{Lemma:smallvolumesimpliessmalldiametersmild}
Let $(M^N,g)$ be a compact Riemannian manifold. There exist two positive constants $\mu^*=\mu^*(M,g)>0$ and $v^*=v^*(M,g)>0$ such that whenever $\Omega\subseteq M$ is an isoperimetric region of volume $0\le v\le v^*$ it holds that $$\diam_g(\Omega)\le\mu^*v^{\frac{1}{N}}.$$  
\end{lemma}

\begin{lemma}\label{Lemma:Concentration+} For any $\eta\in\left]0,1\right[$ sufficiently small and any $r\in\left]0,\frac12{\inj_M}\right[$ there exists $V_2=V_2(M^N, g,\eta, r)>0$ s.t.\ for all $V\in\left]0,V_2\right]$ there exists $\delta_0=\delta_0(\eta,V, M^N, g, W)>0$ such that for every $\delta\in\left]0,\delta_0\right[$ there exists $\varepsilon_2=\varepsilon_2(g, W, \eta, V, \delta)>0$ such that for every 
$\varepsilon\in\left]0,\varepsilon_2\right[$ and for any $u\in \mathfrak M^V_{\varepsilon,c}$ with $c=c(W, V,\delta)=\sigma I_M(V)+\delta$ there exists $\Omega_{V,u}$ isoperimetric region of volume $V$ such that $||u-\chi_{\Omega_{V,u}}||_{L^1(M)}\le\eta$.  In particular for any $\widetilde\eta\in\left]0,1\right[$ $($close to $1$$)$  there exists $p_u\in M$ such that 
\begin{equation}\label{Eq:ConcentrationLowerBound1}
\int_{B_g(p_u,r/2)}|u|dv_g\ge \int_{B_g(p_u,r/2)}udv_g \ge \widetilde\eta V.
\end{equation}
\end{lemma}
\begin{proof} By \cite[Lemma~4.9]{NardulliOsorioIMRN} reported above as Lemma \ref{Lemma:smallvolumesimpliessmalldiametersmild} we know that there exists $v^*_0:=v^*_0(N,k,\inj_M,r)>0$ such that for every isoperimetric region $\Omega$ of volume $V$ smaller than $v^*_0$ is contained in a geodesic ball of radius $r/2$. Furthermore, by Lemma \ref{Lemma:ConcentrationIsoperimetric} we get the existence of a isoperimetric region $\Omega_{V,u}$ such that
$$||u-\chi_{\Omega_{V,u}}||_{L^1(M)}\le\eta.$$
This implies
$$\int_{M\setminus B_g(p_u,r/2)}udv_g\le\int_{M\setminus B_g(p_u,r/2)}|u|dv_g\le\eta,$$
furthermore
$$V=\int_{M\setminus B_g(p_u,r/2)}udv_g+\int_{B_g(p_u,r/2)}udv_g,$$
from which we deduce 
$$\int_{B_g(p_u,r/2)}udv_g\ge V-\eta.$$
From this it is straightforward to choose a suitable $V_2:=V_2(M^N, g,\eta, r)\le v_0^*$ and to deduce \eqref{Eq:ConcentrationLowerBound1} for small $\eta$, and thus we finish the proof of the lemma.
\end{proof}
\begin{remark}
Notice that in Lemma \ref{Lemma:Concentration+} 
if $\eta$ is  small enough we can ensure that  $\frac12 \le\frac{ \int_M  u(x)dv_g}{\int_M|u(x)|dv_g}\le 1$. With this observation it is easy to conclude that in our results remain true when we replace the barycenter  $\beta^*$ by other notions of barycenter as listed bellow
$$
\begin{aligned}
 \beta^*_1(u)=\frac{ \int_M  xu(x)dv_g}{\int_M|u(x)|dv_g},&\qquad &  \beta^*_2(u)=\frac{ \int_M  x|u(x)|dv_g}{\int_M|u(x)|dv_g},\\
 \beta^*_3(u)=\frac{ \int_M  x|u(x)|dv_g}{\int_Mu(x)dv_g}.&\qquad  & 
\end{aligned}
$$
\end{remark}
Let us denote by $\diam_{\R^l}(M)$ the diameter of $M$ as subset of $\R^l$.

\begin{lemma}\label{Lemma:WellPosedBarycenters} For $r\in\left]0,\frac12\inj_M\right[$ there exists $V_4=V_4(N,g,v_0^*, \inj_M, r, \diam_{\R^l}(M))>0$  such that for every $0<V<V_4$, there exists $\varepsilon_4=\varepsilon_4(g, W, V)>0$, $0<\varepsilon<\varepsilon_4$, and every $u\in\mathfrak{M}^V_{\varepsilon,c}$ we have $\beta^*(u)\in M_{r}$. 
\end{lemma}
\begin{proof} Define $\rho(u(x)):=\frac{u(x)}{\int_Mu(x)dv_g}$. By 
\eqref{Eq:ConcentrationLowerBound1} for every $V\in\left]0,V_3\right[$ we obtain a point $p_u\in M$ such that  $\int_{B_g(p_u,r/2)}\rho(u(x))dv_g\ge\tilde{\eta}$, where $0<\tilde{\eta}<1$ will be chosen later. From this last inequality we deduce
\begin{eqnarray*}
|\beta^*(u)-p_u| & = & \left|\int_M(x-p_u)\rho(u(x))dv_g\right|\\
& \le & \left|\int_{B_g(p_u,r/2)}(x-p_u)\rho(u(x))dv_g\right|\\
& + & \left|\int_{M\setminus B_g(p_u,r/2)}(x-p_u)\rho(u(x))dv_g\right|\\
& \le & \frac{r}{2}+D(1-\tilde{\eta}),
\end{eqnarray*}
where $D:=\diam_{\R^l}(M)$.
Choosing $\tilde{\eta}$ close to $1$ such that $D(1-\tilde{\eta})<\frac{r}{2}$ and applying Lemma \ref{Lemma:Concentration+}, with $\eta=\tilde{\eta}$, and $\delta\in]0,\delta_0[$ we conclude the proof of the lemma. 
\end{proof}
\begin{corollary}\label{Cor:HomotopyOnM} There exists $r_0=r_0(M)>0$ such that for any $r\in\left]0,r_0\right[$, there exists $V_5=V_5\big(N,k,v_0^*, \inj_M, r, \diam_{\R^l}(M)\big)>0$ such that for every $V\in\left]0,V_5\right[$, there exists $\varepsilon_5=\varepsilon_5(V)>0$ such that for every $\varepsilon\in\left]0,\varepsilon_5\right[$, we have $d_g(\pi\circ\beta^*\circ\Phi_{\varepsilon,V}(x_0),x_0)<\inj_M$. In particular $\pi\circ\beta^*\circ\Phi_{\varepsilon,V}$ is homotopic to the identity map of $M^N$.
\end{corollary}
\begin{proof} As in Lemma~\ref{Lemma:WellPosedBarycenters}, if we choose $r_0$ small enough depending only on the second fundamental form of the isometric immersion of $M$ in $\R^l$ and the injectivity radius of $M$, it is easy to see that we have $d_g(\pi\circ\beta^*\circ\Phi_{\varepsilon,V}(x_0),x_0)\le C(||II_M||_{\infty})r_0<\inj_M$, because $M$ is compact. To understand this standard argument of extrinsic Riemannian geometry, the reader can look up \cite[Lemma~2.1]{NardulliBBMS2018}. Let us now define the homotopy $F:[0,1]\times M\to M$, 
\[F(t,x_0):=\exp_{x_0}(t\exp^{-1}_{x_0}(\pi\circ\beta^*(\Phi_{\varepsilon,V}(x_0)))).\] 
From the very definition of $F$ it is easy to check that $F(0,x_0)=x_0$ and $F(1,x_0)=\pi\circ\beta^*\circ\Phi_{\varepsilon,V}(x_0)$ for every $x_0\in M$. Checking the continuity of $F$ with respect to $x_0$ is a standard fact of Riemannian geometry about the exponential map using Proposition \ref{Lemma:EnergyUpperBounds} (continuity of $\Phi_{\varepsilon, V}$ ) and Lemma \ref{Lemma:BarycenterContinuity} (continuity of barycenter map).  
\end{proof}

We are finally in position to prove Theorem \ref{Thm:Main1}.

\begin{proof} Set $V^*:=\min\{V_2, V_3, V_4, V_5\}>0$, then fix $0<\delta<\delta_0$, with $\delta_0$ as in Lemma \ref{Lemma:Concentration+}, and set  $\varepsilon^*:=\min\{\varepsilon_0,\varepsilon_1,\varepsilon_2, \varepsilon_3,\varepsilon_4,\varepsilon_5\}>0$, $c=\sigma I_M(V)+\delta$. Then for any $V\in\left]0,V^*\right[$ and $\varepsilon\in\left]0,\varepsilon^*\right[$,  by an easy application of Proposition~\ref{Lemma:EnergyUpperBounds}, Lemma~\ref{Lemma:WellPosedBarycenters} and Corollary~\ref{Cor:HomotopyOnM} 
we obtain the functions $f:=\Phi_{\varepsilon,V}$ and $g:=\pi\circ\beta^*$ required to apply Theorem \ref{Thm:FotoLS} to $X=M$, $J={E_\varepsilon}\vert_{{\mathfrak{M}^V}}$,  $\mathfrak{M}=\mathfrak{M}^V$. The conclusion then follows readily. 

The last assertion of the theorem follows directly from Theorem \ref{Thm:FotoMorse}, using the nondegeneracy assumption.
\end{proof}

\section{Dropping the subcritical growth condition}\label{sec:dropping}
We will now show how to deal with the case where one does not assume the subcritical growth of the potential \eqref{Eq:Potential0}. The idea is to show some \emph{a priori} estimates on the solutions (and for the corresponding Lagrange multiplier), and then consider a perturbed problem that satisfies the growth condition, whose solutions are also solutions of the original problem.
Towards this goal,  we need two auxiliary lemmas that have their own interest. Recalling Remark~\ref{thm:remsuperjacent}, a careful inspection of the proof of Theorem \ref{Thm:Main1},  
Corollary \ref{Cor:HomotopyOnM}, Lemma \ref{Lemma:WellPosedBarycenters} reveals that 
\begin{lemma}\label{Lemma:EnergyBounds} Let $W$ satisfying assumptions \eqref{Eq:Potential}, \eqref{Eq:Potential0}, and \eqref{Eq:WeakUpperBarriers}. Then, for every $V\in]0,V^*[$, $\varepsilon\in]0,\varepsilon^*[$, there exists $\hat c=\hat c\big(N,g, \varepsilon,V,s_0,W_{\vert_{[0,s_0]}}\big)>\inf\nolimits_{\mathfrak M^V} E_\varepsilon$, such that for all $c\in\left]\inf\limits\nolimits_{\mathfrak M^V} E_\varepsilon,\hat c\right]$, the sublevel $E_\varepsilon^c$ is \emph{homotopically superjacent} to $M$ $($see Remark~\ref{thm:remsuperjacent}$)$.\qed
\end{lemma}
The proof of the following result goes along the same lines as the proof of \cite[Lemma~3.4]{Chen96}. Before giving its statement, it will be useful to make a remark.
\begin{remark} Let us observe that from the quadratic growth condition obtained integrating two times $W''(t)\ge c_0>0, \forall |t|\ge t_0$, on the interval $[t_0,s]$ we obtain $W(s)\ge W(t_0)+(s-t_0)W'(t_0)+\frac12(s-t_0)^2c_0$ from which we conclude that there exists $t_1\ge t_0$, $c_0'>c_0$, such that $W(s)\ge\frac12s^2c'_0,$ for every $s$ such that $|s|\ge t_1$.
\end{remark}
\begin{proposition}[Lagrange Multiplier Estimates]\label{Lemma:Lagrange Multiplier Estimates1} Let $\mathcal{E}_0, V$, and $\overline\varepsilon\in]0, +\infty[$ be fixed, and assume that $(u_\varepsilon)_{\varepsilon\in\left]0,\overline\varepsilon\right]}$ is a family of solutions for the equation 
\begin{equation}\label{Eq:UniformLangrageMultiplierEstimates01}
-\dive\left(\varepsilon\nabla u_{\varepsilon}\right)+\tfrac1{\varepsilon}W'\big(u_{\varepsilon}\big)=-\lambda_{\varepsilon},\quad\text{in }M, \int_M u_\varepsilon=V>0,
\end{equation} 
where $W$ satisfy \eqref{Eq:Potential}, $W''(1)>0$, 
$W''(s)\ge c_0>0,\text { if }|u| \geq t_0$ for some $c_0>0$ and $t_0>s_0>0$, i.e., large quadratic or super quadratic growth such that 
\begin{equation}\label{Eq:UniformLangrageMultiplierEstimates02}
0\le E_\varepsilon[u_\varepsilon]\le\mathcal{E}_0, \forall\varepsilon\in\left]0,\overline\varepsilon\right].
\end{equation}
Then there exist positive constants $c_1=c_1(N,g, \mathrm{vol}_g(M), V, \mathcal{E}_0, t_0, c_0, W_{\vert_{[-t_0, t_0]}})>0$ $($$c_1>0$ large$)$ and $\varepsilon_0=\varepsilon_0(N,g, \mathrm{vol}_g(M), V, \mathcal{E}_0, t_0, W_{\vert_{[-t_0, t_0]}})>0$, $($$\varepsilon_0>0$ small$)$ such that for any $\hat{\varepsilon}\in]0,\varepsilon_0]$ we have $c_1\le \widetilde{c}_1$, $\widetilde{c}_1=\widetilde{c}_1(N, g, \mathrm{vol}_g(M), V, \mathcal{E}_0, t_0, c_0, W_{\vert_{[-t_0, t_0]}})>0$ with 
\begin{eqnarray*} 
|\lambda_{\hat{\varepsilon}}(u_{\hat{\varepsilon}})| & \le & c_1E_{\hat{\varepsilon}}(u_{\hat{\varepsilon}})\\
& \le & \widetilde{c}_1E_{\hat{\varepsilon}}(u_{\hat{\varepsilon}})\\
& \le & \widetilde{c}_1\mathcal{E}_0.
\end{eqnarray*}   
\end{proposition}
\begin{remark} The assumptions of Proposition~\ref{Lemma:Lagrange Multiplier Estimates1} are satisfied in the case of the classical symmetric Van der Waals-Allen-Cahn-Hilliard potential that is a positive polynomial of fourth order with just two absolute minima at which the potential is zero. 
\end{remark}
\begin{remark}\label{Rem:LagrangeMultiplierEstimates1ConstantsDependence} 
Roughly speaking, Proposition~\ref{Lemma:Lagrange Multiplier Estimates1} says that the constants involved in the statement of our results depend on the geometry of the problem, on an upper bound of the energy, on the  behavior of the potential over a compact interval, and on the index of quadratic and super-quadratic growth at infinity, which is represented by the constant $c_0>0$. 
\end{remark}
For the needs of the proof of Proposition~\ref{Lemma:Lagrange Multiplier Estimates1} we need the following, by now, standard notion of mollification kernel and of mollification of a function defined on a Riemannian manifold. Along the past decades this central topic has been treated by many authors. At the best of our knowledge the first appearance in the literature, of an intrinsic treatment of mollifier kernels on a smooth Riemannian manifold, dates back to the work of H. Karcher \cite{Karcher77} while for an extrinsic treatment the earliest reference that we found is the celebrated work of John Nash \cite{Nash56Imbedding}. The treatment presented here is borrowed from \cite{Fukuoka06}.
Let $U \subset \subset M$ be an open set (this includes the case $U=M$ if $M$ is a closed differentiable manifold). Define the injectivity radius of $U$ by
$$
\operatorname{inj}(U, g)=\inf _{x \in U} \operatorname{inj}_g(x),
$$
where $\operatorname{inj}_g(x)$ is the injectivity radius of $x \in M$ with respect to the metric $g$. 
\begin{definition} We say that a function $\phi: U \times M \times(0, \operatorname{inj}(U, g))$ is a \textbf{\emph{mollifier kernel}}, if
\begin{enumerate}
\item [(1)] $\phi\in C^\infty(U \times M \times(0, \operatorname{inj}(U, g)))$.
\item [(2)] $\phi(x, \cdot,\rho): M \rightarrow \mathbb{R}$ has its support in $\bar{B}_g(x, \rho)$.
\item [(3)] $\int_{M} \phi(x, y, \rho)dv_g(y)=1$ for every $x \in U$ and $\rho\in(0, \operatorname{inj}(U, g))$ fixed.
\end{enumerate}
\end{definition}
Define
$$
\check{\psi}(x, y, \rho)=\left\{\begin{array}{l}
e^{\frac{1}{\left(\frac{\operatorname{dist}(x, y)}{\rho}\right)^{2}-1}}, \text { if } \operatorname{dist}(x, y)<\rho<\operatorname{inj}(U, g), \\
0, \text { if } \operatorname{dist}(x, y) \geq\rho.
\end{array}\right.
$$

The function $\check{\psi}$ belongs clearly to $C^{\infty}(M)$. Now, let us set
$$
\psi(x, y, \rho):=\frac{\check{\psi}(x, y, \rho)}{\int_{M} \check{\psi}(x, y, \rho) dv_g(y)}, \psi_\rho(x, y):=\psi(x, y, \rho).
$$
\begin{remark}
Notice that $\psi(x, \cdot, \rho)$ has its support in $B_g(x, \rho)$.
\end{remark}
It is not difficult to see that the function $\psi$ defined above satisfies the properties of a mollifier on $(M, g)$ and it will be called the \textbf{\emph{standard mollifier}}. 

For every $f\in L^1_{loc}(M)$ the \textbf{\emph{mollified function}} $f_{\rho}$ is defined as follow
$$f_{\rho}(x):=(f\ast\psi_\rho)(x):=\int_M\psi(x, y, \rho)f(y)dv_g(y).$$ It is easy to see that $f_\rho\in C^\infty(M)$.


Some steps in the proof of Proposition \ref{Lemma:Lagrange Multiplier Estimates1} rely on basic properties of this mollification construction performed on Riemannian manifolds (such as [GT01, Lemma 7.23]) – concretely, the estimates for $C^1$ and $L^\infty$ norms of $u_{\varepsilon,\rho}$, and the upper bound for $||u_{\varepsilon,\rho}-u_{\varepsilon}||_{L^2(M)}$, where $u_{\varepsilon,\rho}$ is a mollification of $u_\varepsilon$ (see inside the body of the proof below for the precise definition of  $u_{\varepsilon,\rho}$). These properties follows mutatis mutandis in a straightforward manner from the corresponding Euclidean ones, for this reason we omit the proof here.
\begin{proof}[Proof of Proposition~\ref{Lemma:Lagrange Multiplier Estimates1}] We can assume w.l.g.\ that $0<\bar{\varepsilon}<1$. Looking at the equation \eqref{Eq:UniformLangrageMultiplierEstimates01} we want to give a uniform estimate with respect to $\varepsilon$ of $\lambda_{\varepsilon,V}$ depending only on the energy of the associated solutions. With this aim in mind, we will make use of an auxiliary function $\psi_{\varepsilon,\rho}:M\to\R$ given as the unique solution to 
 \begin{equation}\label{Eq:UniformLangrageMultiplierEstimatesAuxiliaryTrickyEllipticProblem}
\left\{\begin{aligned}
\Delta\psi_{\varepsilon,\rho} & = u_{\varepsilon,\rho}-\bar{u}_{\varepsilon,\rho}, \text{ in }M,\\
\int_M\psi_{\varepsilon,\rho}dv_g & = 0,
\end{aligned}
\right.
\end{equation}
with $u_{\varepsilon,\rho}(x):=\left(u_\varepsilon\ast\psi_\rho\right)(x):=\int_M\psi(x, y, \rho)u_\varepsilon(y)dv_g(y)$, where $\psi_\rho(x,y):=\psi(x,y,\rho)$ is the standard mollification kernel satisfying $\int\psi_\rho=1$ defined above and $$\bar{u}_{\varepsilon,\rho}:={\mathrm{vol}_g(M)}^{-1}\int_M u_{\varepsilon,\rho}dv_g=\frac V{\mathrm{vol}_g(M)}.$$ By a direct computation coming from the very definition of $u_{\varepsilon,\rho}$ we get
\begin{eqnarray*}
||u_{\varepsilon,\rho}||_{\infty,M} & = & ||[(u_\varepsilon-1)+1]\ast\psi_\rho||_{\infty,M}\\
 & \le & 1+\sup_{x\in M}\int_{B_g(x,\rho)}\psi_\rho(x,y)\left||u_\varepsilon(\exp_{x,g}(-\exp_{x,g}^{-1}(y))|-1\right|\,\mathrm dv_g(y)\\
 & \stackrel{\text{H\"older}}{\le} & 1+C\rho^{-\frac{N}{2}}|||u_\varepsilon|-1||_{2,M}\\
 & \le & 1+C\sqrt{\mathcal{E}_0}\varepsilon^{\frac12}\rho^{-\frac{N}{2}},
\end{eqnarray*}
where the last constant $C=C(N,g)$.
The last inequality is an immediate consequence of 
\begin{equation}\label{Eq:LagrangeMultiplierStartingEstimates}
\int_{M}\left(\left|u_{\varepsilon}\right|-1\right)^{2}dv_g\le C\varepsilon\mathcal{E}_0,
\end{equation} 
where $C=C\big(g,W_{\vert_{[-t_0,t_0]}}, c_0\big)>0$.

In order to show \eqref{Eq:LagrangeMultiplierStartingEstimates} we start by considering the Taylor expansion of $W$ near the point $s=1$ on the $s$-axis, which gives the existence of $\theta$ between $1$ and $s$ such that 
\[ W(s)=W(s)-W(1)=W'(1)(s-1)+\frac12(s-1)^2W''(\theta_s)=\frac12(s-1)^2W''(\theta_s).\]
From this we easily obtain that for $|s-1|\le\delta\le\frac12$ 
\[W''(\theta_s)\ge\eta_0(\delta, W_{\vert_{[1-\delta, 1+\delta]}})>0,\]
because by assumption \eqref{Eq:Potential} we have $W''(1)>0$ and $W$ is of class $C^2$. For $D:=\{s\in\R||s-1|\ge\delta, s\in[-t_0,t_0]\}$ we have that the function $g$ defined as
$$g(s):=W''(\theta_s)=\frac{2W(s)}{(s-1)^2}>0, g:D\to\R,$$ is continuous and stays away from zero on a compact interval.
Hence $\inf_{s\in D}\{g(s)\}=:\eta_1(\delta, W\vert_{D})>0$. Finally for $|s|\ge t_0$ it is immediate to get $W''(\theta_s)>c_0>0$. This argument implies readily that there exists $\eta_2=\eta_2(W_{\vert_{[-t_0,t_0]}}, c_0):=\min\{\eta_0,\eta_1, c_0\}>0$ such that
$$\int_M(|u_\varepsilon|-1)^2\le\int_M(u_\varepsilon-1)^2\le\frac1{\eta_2}[\varepsilon E_\varepsilon[u_\varepsilon]-\frac\varepsilon2||\nabla u_\varepsilon||^2_{2}]\le\frac1{\eta_2}\varepsilon\mathcal{E}_0.$$ From the last inequality we infer quickly \eqref{Eq:LagrangeMultiplierStartingEstimates} setting $C:=\frac1{\eta_2}$.
Analogously it is not too hard to see that from the very definition of the mollifier and the theorem of derivation under the integral sign we get
\begin{equation}\label{Eq:LagrangeMultiplierStartingEstimates03}
||u_{\varepsilon,\rho}||_{C^1(M)}\le C(\nabla_g\psi, \mathrm{vol}_g(M), W_{\vert_{[-t_0,t_0]}}, c_0)\widetilde{\mathcal{E}}_0\rho^{-1}(1+\varepsilon^{\frac12}\rho^{-\frac{N}{2}}),
\end{equation}
where $\widetilde{\mathcal{E}}_0:=\max\{\sqrt{\mathcal{E}_0},1\}$. 
Thus by classical Schauder's elliptic estimates, applied to \eqref{Eq:UniformLangrageMultiplierEstimatesAuxiliaryTrickyEllipticProblem}, taking $\rho$ small enough we conclude 
\begin{eqnarray}\label{Eq:SchauderEstimates}
||\psi_{\varepsilon,\rho}||_{C^2(M)} 
&\le&  C||u_{\varepsilon,\rho}-\bar{u}_{\varepsilon,\rho}||_{C^1(M)}\\ \nonumber
& \le & C\left(||u_{\varepsilon,\rho}||_{C^1(M)}+||\bar{u}_{\varepsilon,\rho}||_{C^1(M)}\right)\\ 
&\stackrel{\eqref{Eq:LagrangeMultiplierStartingEstimates03}}{\le}& 
C(g,\nabla_g\psi, V, \mathrm{vol}_g(M), W_{\vert[-t_0,t_0]}, c_0)\widetilde{\mathcal{E}}_0\rho^{-1}(1+\varepsilon^{\frac12}\rho^{-\frac{N}{2}}).\notag
\end{eqnarray}
Now we come back to our uniform estimates on $\lambda_\varepsilon$ and multiply \eqref{Eq:UniformLangrageMultiplierEstimates01}, by the function $\varphi_{\varepsilon,\rho}:=\langle\nabla\psi_{\varepsilon,\rho},\nabla u_\varepsilon\rangle_g$ then we integrate over $M$ and use the divergence theorem obtaining
\begin{eqnarray}\label{Eq:ChenEstimate0}
\begin{aligned}
\int_M\varphi_{\varepsilon,\rho}(-\lambda_\varepsilon)dv_g
= & \int_M\varphi_{\varepsilon,\rho}\left(-\dive(\varepsilon\nabla u_\varepsilon)+\frac1{\varepsilon}W'(u_\varepsilon)\right)dv_g\\
= & \varepsilon \int_{M}\left\langle\nabla \varphi_{\varepsilon, \rho}, \nabla u_{\varepsilon}\right\rangle_gdv_g+\frac1{\varepsilon}\int_{M}W'(u_\varepsilon)\langle\nabla\psi_{\varepsilon,\rho},\nabla u_\varepsilon\rangle_gdv_g\\
= &\int_{M}\frac\varepsilon2|\nabla u_\varepsilon|^2\dive(\nabla\psi_{\varepsilon,\rho})dv_g -\frac1{\varepsilon}\int_{M}W(u_\varepsilon)\dive(\nabla\psi_{\varepsilon,\rho})dv_g\\ 
\le & \int_{M}\left(\frac\varepsilon2|\nabla u_\varepsilon|^2+\frac1{\varepsilon}W(u_\varepsilon)\right)|\dive(\nabla\psi_{\varepsilon,\rho})|dv_g\\
\le & ||\psi_{\varepsilon,\rho}||_{C^2(M)}E_{\varepsilon}(u_\varepsilon)\\ \label{Eq:SchauderEstimates0}
\stackrel{\eqref{Eq:SchauderEstimates}}{\le} & C\widetilde{\mathcal{E}}_0\rho^{-1}(1+\varepsilon^{\frac12}\rho^{-\frac{N}{2}})E_{\varepsilon}(u_\varepsilon)
\end{aligned}
\end{eqnarray}
where $C(g,\nabla_g\psi, \mathrm{vol}_g(M), W_{\vert[-t_0,t_0]}, c_0)>0$.  We will explain in great details in the next paragraph the identity displayed in the third line of \eqref{Eq:ChenEstimate0}. With this aim in mind observe that an integration by parts gives 
\begin{eqnarray}\label{Eq:Identity0}
\int_M \langle\nabla\varphi_{\varepsilon,\rho},\nabla u_\varepsilon\rangle dv_g & = & -\frac12\int_Mu_\varepsilon\Delta\varphi_{\varepsilon,\rho}dv_g-\frac12\int_M\varphi_{\varepsilon,\rho}\Delta u_\varepsilon dv_g.
\end{eqnarray}
Remembering of the following standard equality 
$$
\nabla\langle X,Y\rangle_g=\dive_g(X)Y+\dive_g(Y)X,
$$
valid for any pair of smooth enough vector fields $X,Y$, and applying it to our specific case with $X=\nabla u_\varepsilon$, $Y=\nabla\psi_{\varepsilon, \rho}$ it holds
\begin{eqnarray*}
\Delta\varphi_{\varepsilon, \rho}=2\Delta\psi_{\varepsilon, \rho}\Delta u_\varepsilon+\langle\nabla\Delta\psi_{\varepsilon, \rho},\nabla u_\varepsilon\rangle_g+\langle\nabla\Delta u_\varepsilon, \nabla\psi_{\varepsilon, \rho}\rangle_g. 
\end{eqnarray*}
Thus
\begin{eqnarray*} 
-\frac12\int_Mu_\varepsilon\Delta\varphi_{\varepsilon, \rho}dv_g & = & -\frac12\int_M\left(2u_\varepsilon\Delta\psi_{\varepsilon, \rho}\Delta u_\varepsilon+\langle\nabla\Delta\psi_{\varepsilon, \rho},u_\varepsilon\nabla u_\varepsilon\rangle_g+\langle\nabla\Delta u_\varepsilon, u_\varepsilon\nabla\psi_{\varepsilon, \rho}\rangle_g\right)dv_g\\
& = & -\frac12\int_M\left(2u_\varepsilon\Delta\psi_{\varepsilon, \rho}\Delta u_\varepsilon-\Delta\psi_{\varepsilon, \rho}\dive (u_\varepsilon\nabla u_\varepsilon)-\Delta u_\varepsilon\dive(u_\varepsilon\nabla\psi_{\varepsilon, \rho})\right)dv_g\\
& = & -\frac12\int_M\left(2u_\varepsilon\Delta\psi_{\varepsilon, \rho}\Delta u_\varepsilon-\Delta\psi_{\varepsilon, \rho}|\nabla u_\varepsilon|^2-2u_\varepsilon\Delta\psi_{\varepsilon, \rho}\Delta u_\varepsilon-\Delta u_\varepsilon\langle\nabla u_\varepsilon,\nabla\psi_{\varepsilon, \rho}\rangle_g\right)dv_g\\
& = & \frac12\int_M\left(\Delta\psi_{\varepsilon, \rho}|\nabla u_\varepsilon|^2+\Delta u_\varepsilon\langle\nabla u_\varepsilon,\nabla\psi_{\varepsilon, \rho}\rangle_g\right)dv_g.
\end{eqnarray*}
Hence substituting in \eqref{Eq:Identity0} we get
\begin{eqnarray*}
\int_M \langle\nabla\varphi_{\varepsilon, \rho},\nabla u_\varepsilon\rangle_gdv_g & = & \frac12\int_M\left(\Delta\psi_{\varepsilon, \rho}|\nabla u_\varepsilon|^2+\Delta u_\varepsilon\langle\nabla u_\varepsilon,\nabla\psi_{\varepsilon, \rho}\rangle_g\right)dv_g-\frac12\int_M\varphi_{\varepsilon, \rho}\Delta u_\varepsilon dv_g\\
& = & \frac12\int_M\left(\Delta\psi_{\varepsilon, \rho}|\nabla u_\varepsilon|^2+\varphi_{\varepsilon, \rho}\Delta u_\varepsilon\right)dv_g-\frac12\int_M\varphi_{\varepsilon, \rho}\Delta u_\varepsilon dv_g\\
& = & \frac12\int_M\Delta\psi_{\varepsilon, \rho}|\nabla u_\varepsilon|^2dv_g.
\end{eqnarray*}
The last identity readily proves the third line of \eqref{Eq:Identity0}.
An integration by parts on the left-hand side of the above inequality \eqref{Eq:ChenEstimate0} yields
\begin{eqnarray*}
\int_M\langle\nabla\psi_{\varepsilon,\rho},\nabla u_\varepsilon\rangle_g dv_g & = & \int_M u_\varepsilon \dive(\nabla\psi_{\varepsilon,\rho})dv_g\\
& = & \int_M u_\varepsilon\Delta\psi_{\varepsilon,\rho}dv_g\\
& \stackrel{\eqref{Eq:UniformLangrageMultiplierEstimatesAuxiliaryTrickyEllipticProblem}}{=} & \int_M u_\varepsilon\left(u_{\varepsilon,\rho}-\bar{u}_{\varepsilon,\rho}\right)dv_g\\
& = & \int_M u_\varepsilon\left(u_{\varepsilon,\rho}-u_\varepsilon\right)dv_g+\int_M (u_\varepsilon^2-1)dv_g\\
&  & +\mathrm{vol}_g(M)(1-\bar{u}_\varepsilon^2)+\mathrm{vol}_g(M)\bar{u}_\varepsilon(\bar{u}_\varepsilon-\bar{u}_{\varepsilon,\rho}). 
\end{eqnarray*}
Recall here that $\bar{u}_\varepsilon=\frac{V}{\mathrm{vol}_g(M)}\in\left]0,1\right[$ when $V\in]0, \mathrm{vol}_g(M)[$; using the equality $(x-1)^2+2(x-1)=x^2-1$ and an application of H\"older inequality we obtain 
\begin{equation}\label{Eq:5.20}
\begin{aligned} 
\int_M|u_\varepsilon^2-1|dv_g & \le  \int_M|u_\varepsilon-1|^2dv_g+2\int_M|u_\varepsilon-1|dv_g\\ 
& \stackrel{\text{\eqref{Eq:LagrangeMultiplierStartingEstimates}}}{\le} C\varepsilon\mathcal{E}_0+2\sqrt{C\varepsilon\mathcal{E}_0}\mathrm{vol}_g(M)^{\frac12}\\ 
& \le  C\hat{\mathcal{E}}_0\sqrt{\varepsilon};
\end{aligned}
\end{equation}
for the last inequality in \eqref{Eq:5.20} we have taken $\hat{\mathcal{E}}_0:=\mathcal{E}_0+2\sqrt{\mathcal{E}_0}\mathrm{vol}_g(M)^{\frac12}>0$, assuming without loss of generality $\varepsilon\in\left]0,1\right[$ and $C>1$. In order to verify
\begin{eqnarray}\label{Eq:5.18}
|\bar{u}_{\varepsilon,\rho}-\bar{u}_\varepsilon| & \stackrel{\text{H\"older}}{\le} &  \mathrm{vol}_g(M)^{-\frac12}||u_{\varepsilon,\rho}-u_\varepsilon||_{2,M}\\\label{Eq:5.18+}
& \le & C\sqrt{\rho},
\end{eqnarray}
where we used \eqref{Eq:L^2EstimatesWithMean} whose proof is independent of that of \eqref{Eq:5.18}, \eqref{Eq:5.18+}.
At this stage, it is convenient to introduce a new auxiliary function $w_{\varepsilon}$ defined by $w_{\varepsilon}=\widetilde{W}\circ u_{\varepsilon}$ where 
\[\widetilde{W}(s)=\int_{0}^{s} \sqrt{2 \widetilde{F}(t)}dt,\]
\[\widetilde{F}(t) :=\min \left\{W(t), 1+|t|^{2}\right\}\ge 1, \quad\forall t \in \R.\]
Notice that $w_\varepsilon\in W^{1,1}(M)$ and 
\[\int_{M}\left|\nabla w_{\varepsilon}\right|dv_g = \int_{M} \sqrt{2 \widetilde{F}\left(u_{\varepsilon}\right) |} \nabla u_{\varepsilon} | dv_g
\le \int_{M}e_{\varepsilon}\left(u_{\varepsilon}\right)dv_g= E_{\varepsilon}[u_\varepsilon]\le \mathcal{E}_0,\]
where $e_\varepsilon(u)=\frac\varepsilon2\vert\nabla u\vert^2+\frac1\varepsilon W(u)$ is the energy density.
Furthermore, by the structural properties of $W$, there is a positive constants \(c^*_{1}>0\) 
, depending only on the structural assumptions on $W$, such that 
\begin{equation}\label{Eq:05}
\left|s_{1}-s_{2}\right|^{2} \leq c^*_{1}\left|\widetilde{W}\left(s_{1}\right)-\widetilde{W}\left(s_{2}\right)\right| 
,\forall s_1,s_2\in\R,
\end{equation}
and
\begin{equation}\label{Eq:L^2EstimatesWithMean}
\tiny
\begin{aligned} 
\int_{M}\left|u_{\varepsilon,\rho}-u_{\varepsilon}\right|^{2} dv_g& \le  \int_{M} \int_{B_{g}(x,\rho)} \psi_\rho(y)\left|u_{\varepsilon}(\exp_{x,g}(-\exp_{x,g}^{-1}(y))-u_{\varepsilon}(x)\right|^{2}\,\mathrm dv_g(y)\,\mathrm dv_g(x)\\
& \stackrel{\eqref{Eq:05}}{\le} c^*_{1}\int_{M} \int_{B_{g}(x,\rho)} \psi_\rho(y)\left|w_{\varepsilon}(\exp_{x,g}(-\exp_{x,g}^{-1}(y))-w_{\varepsilon}(x)\right|\,\mathrm dv_g(y)\,\mathrm dv_g(x)\\ 
& \stackrel{Fubini}{\le} c^*_{1}\int_{B_{g}(x,\rho)}\psi_\rho(y)\int_{M} \left|w_{\varepsilon}(\exp_{x,g}(-\exp_{x,g}^{-1}(y))-w_{\varepsilon}(x)\right|\,\mathrm dv_g(x)\,\mathrm dv_g(y)\\ 
& \stackrel{\text{Riemannian version of \cite[Lemma~7.23]{GilbargTrudinger98}}}{\le}  c^*_1\rho\left\|\nabla w_{\varepsilon}(\cdot)\right\|_{1, M} \leq C \rho,
\end{aligned}
\normalsize
\end{equation}
where $C=C(g,\mathcal{E}_0, \mathrm{vol}_g(M))>0$. From the last inequality and \eqref{Eq:5.18} we obtain easily \eqref{Eq:5.18+}.
Note that $$\int_{M} u_{\varepsilon}^{2} \leq \int_{M}\left|u_{\varepsilon}^{2}-1\right|+\mathrm{vol}_{g}(M),$$ and then we use \eqref{Eq:5.20} to give the following estimate on the $L^2$ norm of $u_\varepsilon$
\begin{equation}
\begin{aligned}\label{Eq:5.19}
\int_M u_\varepsilon^2 & \le C\hat{\mathcal{E}}_0\sqrt{\varepsilon}+\mathrm{vol}_{g}(M),
\end{aligned} 
\end{equation}
which implies by an application of H\"older inequality that 
\small
\begin{eqnarray}\nonumber
\left|\int_M u_\varepsilon\left(u_{\varepsilon,\rho}-u_\varepsilon\right)dx\right| & \le & ||u_\varepsilon||_2 ||u_{\varepsilon,\rho}-u_\varepsilon||_2\\ \label{Eq:5.22}
 & \stackrel{\eqref{Eq:5.18}-\eqref{Eq:5.19}}{\le} & C\sqrt{C\hat{\mathcal{E}}_0\sqrt{\varepsilon}+\mathrm{vol}_{g}(M)}\sqrt{\rho}=C\mathcal{E}_0^*\sqrt{\rho}.
\end{eqnarray}
\normalsize
So 
\begin{multline}\label{Eq:AuxiliaryProblemSchauder}
\int_M u_\varepsilon \dive(\nabla\psi_{\varepsilon,\rho})dx\\
\stackrel{\eqref{Eq:5.20}-\eqref{Eq:5.22}}{\ge} \mathrm{vol}_g(M)\left(1-\left(\frac{V}{\mathrm{vol}_g(M)}\right)^2\right)-C\hat{\mathcal{E}}_0(\sqrt{\varepsilon})-C(1+\mathcal{E}_0^*)\sqrt{\rho}.
\end{multline}
Now combining \eqref{Eq:SchauderEstimates0} and \eqref{Eq:AuxiliaryProblemSchauder}, for sufficiently small $\varepsilon$ and $\rho$, depending only on the relevant quantities $g, N, \operatorname{vol}_{g}(M), V, \mathcal{E}_{0}, t_{0}$, $c_{0}$, and $W_{\vert_{\left[-t_0, t_0\right]}}$, we deduce that 
\begin{equation}
|\lambda_{\varepsilon,V}|\le\frac{C\widetilde{\mathcal{E}}_0\rho^{-1}(1+\varepsilon^{\frac12}\rho^{-\frac{N}{2}})E_{\varepsilon}(u_\varepsilon)}{\mathrm{vol}_g(M)(1-\left(\frac{V}{\mathrm{vol}_g(M)}\right)^2)-C\hat{\mathcal{E}}_0\sqrt{\varepsilon}-C(1+\mathcal{E}_0^*)\sqrt{\rho}}.
\end{equation}
So taking $\varepsilon$ such that 
\begin{equation}
\frac12\mathrm{vol}_g(M)\left(1-\left(\frac{V}{\mathrm{vol}_g(M)}\right)^2\right)-C\hat{\mathcal{E}}_0\sqrt{\varepsilon}\le \frac{\mathrm{vol}_g(M)}{4}\left(1-\left(\frac{V}{\mathrm{vol}_g(M)}\right)^2\right),
\end{equation}
and $\rho$ such that
\[
\frac12\mathrm{vol}_g(M)\left(1-\left(\frac{V}{\mathrm{vol}_g(M)}\right)^2\right)-C(1+\mathcal{E}_0^*)\sqrt{\rho}\le\frac{\mathrm{vol}_g(M)}{4}\left(1-\left(\frac{V}{\mathrm{vol}_g(M)}\right)^2\right),
\]
we conclude that 
\[
|\lambda_{\varepsilon,V}|\le\frac{2C\widetilde{\mathcal{E}}_0\rho^{-1}(1+\varepsilon^{\frac12}\rho^{-\frac{N}{2}})E_{\varepsilon}(u_\varepsilon)}{\mathrm{vol}_g(M)\left(1-\left(\frac{V}{\mathrm{vol}_g(M)}\right)^2\right)}.\qedhere
\]
\end{proof}

\begin{remark} In Proposition~\ref{Lemma:Lagrange Multiplier Estimates1} we do not require any subcritical growth condition, just that the second derivative of the potential is bounded from below by a positive constant in a neighborhood of infinity.
\end{remark}
In our next result we generalize Theorem \ref{Thm:Main1} by dropping assumption \eqref{Eq:Potential0}, and replacing it with a quite more general one. The price we pay is a weaker estimate on the lower bound on the number of solutions than the one determined in Theorem~\ref{Thm:Main1}. In fact, our proof gives only the existence of low energy solutions. We conjecture that the following result  remains true also for high energy solutions, but at the moment we are unable to give a complete proof. Thus under the assumptions of the preceding lemma we state the following theorem.
\begin{theorem}\label{Thm:Main1Generalized+}
For every $W\in C^2(\mathds  R)$ satisfying \eqref{Eq:Potential} and $W''(1)>0$, $W''(s)>c_0>0,\forall |s|\ge t_1>1>0$, for some large $t_1$, then for every $V\in\left]0,V^*\right[$ and $\varepsilon\in\left]0,\varepsilon^*\right[$, Problem \eqref{Eq:EulerLagrangeVEpsilon} admits at least $\cat(M)$ distinct solutions.
Moreover, assume additionally that for given $V\in\left]0,V^*\right[$ and $\varepsilon\in\left]0,\varepsilon^*\right[$ all solutions of  Problem \eqref{Eq:EulerLagrangeVEpsilon} having energy less than or equal to the constant $\hat{c}$ defined in Lemma \ref{Lemma:EnergyBounds} are nondegenerate $($see Definition~\ref{Eq:DefTopologicallyNondegenerate}$)$. Then,   Problem \eqref{Eq:EulerLagrangeVEpsilon} has at least $P_{1}(M)$ distinct solutions.
\end{theorem}
\begin{proof} Let $V^*$ be the constant determined in Theorem~\ref{Thm:Main1}, and assume in the rest of the proof that $V\in\left]0,V^*\right[$.
We can suppose that $W$ satisfies 
\begin{equation}\label{Eq:FinalPotentialGrowth}
\limsup_{s\to+\infty} W'(s)=+\infty, \text{ and }\liminf_{s\to-\infty}W'(s)=-\infty,
\end{equation} 
because otherwise $W'$ would be bounded, and so $W$ would satisfy a growth condition as in \eqref{Eq:Potential0}, falling under the assumptions of Theorem \ref{Thm:Main1}. Now, using \eqref{Eq:FinalPotentialGrowth} it is easy to see that there exists $\hat{s}^-\le-t_1$, $\hat{s}^+\ge t_1$ such that 
\begin{equation}\label{Eq:SupremumNormBoundsOnSolutions}
-\tfrac1\varepsilon\, W'(\hat{s}^-)-\lambda^*>0,
\end{equation}
and
\begin{equation}\label{Eq:SupremumNormBoundsOnSolutions0}
-\tfrac1\varepsilon\, W'(\hat{s}^+)+\lambda^*<0,
\end{equation} 
where $0<V<V^*$, $\lambda^*=\lambda^*(N, g, \mathrm{vol}_g(M), \varepsilon, V, s_0, t_0, \hat{c}, W_{\vert_{[-t_1,t_1]}}):=c_1\hat{c}>0$ with the notations of Lemma ~\ref{Lemma:EnergyBounds} and Proposition~\ref{Lemma:Lagrange Multiplier Estimates1}. Consider the quadratic truncated problem \eqref{Eq:ConstrainedAllenCahnTruncated}: for fixed positive constants $V$ and $\varepsilon$, find $u\in H^{1}(M)$, and $\lambda \in \mathds  R$ such that
\begin{equation}
\begin{aligned}\label{Eq:ConstrainedAllenCahnTruncated}
-&\varepsilon ^{2}\Delta u+\widehat{W}^{\prime }(u) =\varepsilon\lambda,\\&
\int_{M }u(x)\,\mathrm dx =V,  
\end{aligned}
\end{equation}
with the same $M$ as in the statement of the theorem and $\widehat{W}\in C^2(\mathds  R)$ satisfying $\widehat{W}(s):=W(s),\:\forall s\in[\hat{s}^-, \hat{s}^+]$,  \eqref{Eq:Potential0}, 
\begin{equation}\label{Eq:GrowthForLowerBarriersAuxiliary}
\phantom{\quad\forall s\in\left]-\infty, \hat{s}^-\right],}
-\tfrac1\varepsilon\widehat W'(s)-\lambda^*>0,\quad\forall s\in\left]-\infty, \hat{s}^-\right],
\end{equation}
\begin{equation}\label{Eq:GrowthForUpperBarriersAuxiliary}
\phantom{\quad\forall s\in\left]-\infty, \hat{s}^-\right],}
-\tfrac1\varepsilon\widehat W'(s)+\lambda^*<0,\quad\forall s\in\left[\hat{s}^+,+\infty\right[.
\end{equation}
 Observe that it is always possible to find such a $\widehat{W}$. It is straightforward to check that Problem \eqref{Eq:EulerLagrangeVEpsilon} satisfies the hypothesis of Theorem \ref{Thm:Main1} and Proposition~\ref{Lemma:Lagrange Multiplier Estimates1}. Furthermore by the very definition of $\widehat W$ we have $\widehat{W}_{\vert_{[\hat{s}^+,\hat{s}^-]}}=W_{\vert_{[\hat{s}^+,\hat{s}^-]}}$. 
We claim that all the solutions with energy (w.r.t.\ $\widehat W$) less than or equal to $\hat{c}$ of Problem \eqref{Eq:ConstrainedAllenCahnTruncated} are also solutions of Problem \eqref{Eq:EulerLagrangeVEpsilon} with energy (w.r.t.\ $W$) less than or equal to $\hat{c}$. Suppose that $(\hat{u}_1, \hat{\lambda}_1)$ is a solution of Problem \eqref{Eq:ConstrainedAllenCahnTruncated} then again standard elliptic regularity theory (compare Theorem $6.19$ of \cite{GilbargTrudinger98}) shows that $\hat{u}_1$ is of class $C^{2,\alpha}_{\text{loc}}(M)$ and using Lemma \ref{Lemma:EnergyBounds}, Proposition~\ref{Lemma:Lagrange Multiplier Estimates1}, inequalities \eqref{Eq:GrowthForLowerBarriersAuxiliary}, and \eqref{Eq:GrowthForUpperBarriersAuxiliary} combined with the maximum principle, it is easy to check that $\hat u_1\in\left[\hat{s}^-
,\hat{s}^+\right]$, so $(\hat{u}_1, \hat{\lambda}_1)$ is also a solution of Problem \eqref{Eq:EulerLagrangeVEpsilon}, since $W$ and $\widehat{W}$ coincide on the interval $[\hat{s}^-, \hat{s}^+]$.  With this last argument, we conclude the proof of the theorem.
\end{proof}


\begin{thebibliography}{99}
	
	\bibitem[ADPM11]{AmbrosioDaPratoMennucci}
	Luigi Ambrosio, Giuseppe Da~Prato, and Andrea Mennucci.
	\newblock {\em Introduction to measure theory and integration}, volume~10 of
	{\em Appunti. Scuola Normale Superiore di Pisa (Nuova Serie) [Lecture Notes.
		Scuola Normale Superiore di Pisa (New Series)]}.
	\newblock Edizioni della Normale, Pisa, 2011.
	
	\bibitem[AFP00]{AmbrosioFuscoPallara}
	L.~Ambrosio, N.~Fusco, and D.~Pallara.
	\newblock {\em Functions of Bounded Variation and Free Discontinuity Problems}.
	\newblock Oxford Science Publications. Clarendon Press, 2000.
	
	\bibitem[Bal90]{Baldo}
	Sisto Baldo.
	\newblock Minimal interface criterion for phase transitions in mixtures of
	{C}ahn-{H}illiard fluids.
	\newblock {\em Ann. Inst. H. Poincar\'e Anal. Non Lin\'eaire}, 7(2):67--90,
	1990.
	
	\bibitem[BBM07]{BenciBonannoMicheletti}
	Vieri Benci, Claudio Bonanno, and Anna~Maria Micheletti.
	\newblock On the multiplicity of solutions of a nonlinear elliptic problem on
	{R}iemannian manifolds.
	\newblock {\em J. Funct. Anal.}, 252(2):464--489, 2007.
	
	\bibitem[BC91]{BenciCerami}
	Vieri Benci and Giovanna Cerami.
	\newblock The effect of the domain topology on the number of positive solutions
	of nonlinear elliptic problems.
	\newblock {\em Arch. Rational Mech. Anal.}, 114(1):79--93, 1991.
	
	\bibitem[BC94]{BenciCeramiCalcVar}
	Vieri Benci and Giovanna Cerami.
	\newblock Multiple positive solutions of some elliptic problems via the {M}orse
	theory and the domain topology.
	\newblock {\em Calc. Var. Partial Differential Equations}, 2(1):29--48, 1994.
	
	\bibitem[BCP91]{BenciCeramiPassaseo}
	V.~Benci, G.~Cerami, and D.~Passaseo.
	\newblock On the number of the positive solutions of some nonlinear elliptic
	problems.
	\newblock In {\em Nonlinear analysis}, Sc. Norm. Super. di Pisa Quaderni, pages
	93--107. Scuola Norm. Sup., Pisa, 1991.
	
	\bibitem[BCW19]{BCW19}
	Costante Bellettini, Otis Chodosh, and Neshan Wickramasekera.
	\newblock Curvature estimates and sheeting theorems for weakly stable cmc
	hypersurfaces.
	\newblock {\em Advances in Mathematics}, 352:133--157, 2019.
	
	\bibitem[Ben95]{BenciNewApproach}
	Vieri Benci.
	\newblock Introduction to {M}orse theory: a new approach.
	\newblock In {\em Topological nonlinear analysis}, volume~15 of {\em Progr.
		Nonlinear Differential Equations Appl.}, pages 37--177. Birkh\"auser Boston,
	Boston, MA, 1995.
	
	\bibitem[BJ61]{BartleJoichi}
	Robert~G. Bartle and James~T. Joichi.
	\newblock The preservation of convergence of measurable functions under
	composition.
	\newblock {\em Proc. Amer. Math. Soc.}, 12:122--126, 1961.
	
	\bibitem[BLO97]{BertiniLandimOlla}
	L.~Bertini, C.~Landim, and S.~Olla.
	\newblock Derivation of {C}ahn-{H}illiard equations from {G}inzburg-{L}andau
	models.
	\newblock {\em J. Statist. Phys.}, 88(1-2):365--381, 1997.
	
	\bibitem[BNP20]{BeNaPiCalcVar2020}
	Vieri Benci, Stefano Nardulli, and Paolo Piccione.
	\newblock Multiple solutions for the van der
	{W}aals--{A}llen--{C}ahn--{H}illiard equation th a volume constraint.
	\newblock {\em Calc. Var. Partial Differential Equations}, 59(2):Paper No. 64,
	30, 2020.
	
	\bibitem[Bog07]{Bogachev2007measure}
	Vladimir~I Bogachev.
	\newblock {\em Measure theory}, volume~1.
	\newblock Springer Science \& Business Media, 2007.
	
	\bibitem[BW18]{BW2018Stable}
	Costante Bellettini and Neshan Wickramasekera.
	\newblock Stable cmc integral varifolds of codimension $1 $: regularity and
	compactness.
	\newblock {\em arXiv preprint arXiv:1802.00377}, 2018.
	
	\bibitem[BW20]{BW20}
	Costante Bellettini and Neshan Wickramasekera.
	\newblock The inhomogeneous allen--cahn equation and the existence of
	prescribed-mean-curvature hypersurfaces.
	\newblock {\em arXiv preprint arXiv:2010.05847}, 2020.
	
	\bibitem[CH58]{CahnHilliard1958}
	John~Werner Cahn and J.~E. Hilliard.
	\newblock Free energy of a non-uniform system i:interfacial energy.
	\newblock {\em J. Chem. Phys.}, 27:258--266, 1958.
	
	\bibitem[Che96]{Chen96}
	Xinfu Chen.
	\newblock Global asymptotic limit of solutions of the {C}ahn-{H}illiard
	equation.
	\newblock {\em J. Differential Geom.}, 44(2):262--311, 1996.
	
	\bibitem[CM18]{ChodoshMantoulidis2018Minimal}
	Otis Chodosh and Christos Mantoulidis.
	\newblock Minimal surfaces and the allen-cahn equation on 3-manifolds: index,
	multiplicity, and curvature estimates.
	\newblock {\em arXiv preprint arXiv:1803.02716}, 2018.
	
	\bibitem[Dey19]{arXiv:1910.00989}
	Akashdeep Dey.
	\newblock Existence of multiple closed cmc hypersurfaces with small mean
	curvature.
	\newblock {\em arXiv:1910.00989 [math.DG]}, 2019.
	
	\bibitem[Dey20]{arXiv:2004.05120}
	Akashdeep Dey.
	\newblock A comparison of the almgren-pitts and the allen-cahn min-max theory.
	\newblock {\em arXiv:2004.05120 [math.DG]}, 2020.
	
	\bibitem[Fed69]{Federer}
	Herbert Federer.
	\newblock {\em Geometric measure theory}.
	\newblock Die Grundlehren der mathematischen Wissenschaften, Band 153.
	Springer-Verlag New York Inc., New York, 1969.
	
	\bibitem[Fuk06]{Fukuoka06}
	Ryuichi Fukuoka.
	\newblock Mollifier smoothing of tensor fields on differentiable manifolds and
	applications to riemannian geometry.
	\newblock {\em arXiv preprint math/0608230}, 2006.
	
	\bibitem[GG18]{MR3814054}
	Pedro Gaspar and Marco A.~M. Guaraco.
	\newblock The {A}llen-{C}ahn equation on closed manifolds.
	\newblock {\em Calc. Var. Partial Differential Equations}, 57(4):Art. 101, 42,
	2018.
	
	\bibitem[GL97]{GiacominLebowitz1}
	Giambattista Giacomin and Joel~L. Lebowitz.
	\newblock Phase segregation dynamics in particle systems th long range
	interactions. {I}. {M}acroscopic limits.
	\newblock {\em J. Statist. Phys.}, 87(1-2):37--61, 1997.
	
	\bibitem[GL98]{GiacominLebowitz2}
	Giambattista Giacomin and Joel~L. Lebowitz.
	\newblock Phase segregation dynamics in particle systems with long range
	interactions. {II}. {I}nterface motion.
	\newblock {\em SIAM J. Appl. Math.}, 58(6):1707--1729, 1998.
	
	\bibitem[GT01]{GilbargTrudinger98}
	David Gilbarg and Neil~S. Trudinger.
	\newblock {\em Elliptic partial differential equations of second order}.
	\newblock Classics in Mathematics. Springer-Verlag, Berlin, 2001.
	\newblock Reprint of the 1998 edition.
	
	\bibitem[Gua18]{MR3743704}
	Marco A.~M. Guaraco.
	\newblock Min-max for phase transitions and the existence of embedded minimal
	hypersurfaces.
	\newblock {\em J. Differential Geom.}, 108(1):91--133, 2018.
	
	\bibitem[HT00]{HutchinsonTonegawa}
	John~E. Hutchinson and Yoshihiro Tonegawa.
	\newblock Convergence of phase interfaces in the van der
	{W}aals-{C}ahn-{H}illiard theory.
	\newblock {\em Calc. Var. Partial Differential Equations}, 10(1):49--84, 2000.
	
	\bibitem[Kar77]{Karcher77}
	Hermann Karcher.
	\newblock Riemannian center of mass and mollifier smoothing.
	\newblock {\em Communications on pure and applied mathematics}, 30(5):509--541,
	1977.
	
	\bibitem[MJ00]{MorganJohnson2000}
	Frank Morgan and David~L. Johnson.
	\newblock Some sharp isoperimetric theorems for {R}iemannian manifolds.
	\newblock {\em Indiana Univ. Math. J.}, 49(3):1017--1041, 2000.
	
	\bibitem[Mod87]{Modica}
	Luciano Modica.
	\newblock The gradient theory of phase transitions and the minimal interface
	criterion.
	\newblock {\em Arch. Rational Mech. Anal.}, 98(2):123--142, 1987.
	
	\bibitem[MPPP07]{MPPP2007}
	M.~Miranda, Jr., D.~Pallara, F.~Paronetto, and M.~Preunkert.
	\newblock Heat semigroup and functions of bounded variation on {R}iemannian
	manifolds.
	\newblock {\em J. Reine Angew. Math.}, 613:99--119, 2007.
	
	\bibitem[Mur09]{MurrayBook}
	J.~D. Murray.
	\newblock {\em Pre-pattern formation mechanism for animal coat markings}.
	\newblock Theoretical and Mathematical Physics. Springer, Berlin, 2009.
	
	\bibitem[Nar09]{NarAnn}
	Stefano Nardulli.
	\newblock The isoperimetric profile of a smooth {R}iemannian manifold for small
	volumes.
	\newblock {\em Ann. Global Anal. Geom.}, 36(2):111--131, 2009.
	
	\bibitem[Nar14]{NarCalcVar}
	Stefano Nardulli.
	\newblock The isoperimetric profile of a noncompact {R}iemannian manifold for
	small volumes.
	\newblock {\em Calc. Var. Partial Differential Equations}, 49(1-2):173--195,
	2014.
	
	\bibitem[Nar18]{NardulliBBMS2018}
	Stefano Nardulli.
	\newblock Regularity of isoperimetric regions that are close to a smooth
	manifold.
	\newblock {\em Bull. Braz. Math. Soc. (N.S.)}, 49(2):199--260, 2018.
	
	\bibitem[Nas56]{Nash56Imbedding}
	John Nash.
	\newblock The imbedding problem for {R}iemannian manifolds.
	\newblock {\em Ann. of Math. (2)}, 63:20--63, 1956.
	
	\bibitem[NOA18]{NardulliOsorioIMRN}
	Stefano Nardulli and Luis~Eduardo Osorio~Acevedo.
	\newblock {Sharp Isoperimetric Inequalities for Small Volumes in Complete
		Noncompact Riemannian Manifolds of Bounded Geometry Involving the Scalar
		Curvature}.
	\newblock {\em International Mathematics Research Notices}, 06 2018.
	
	\bibitem[Pit81]{MR626027}
	Jon~T. Pitts.
	\newblock {\em Existence and regularity of minimal surfaces on {R}iemannian
		manifolds}, volume~27 of {\em Mathematical Notes}.
	\newblock Princeton University Press, Princeton, N.J.; University of Tokyo
	Press, Tokyo, 1981.
	
	\bibitem[PP10]{Pisante2010}
	Adriano Pisante and Marcello Ponsiglione.
	\newblock Phase transitions and minimal hypersurfaces in hyperbolic space.
	\newblock {\em Communications in Partial Differential Equations},
	36(5):819--849, 2010.
	
	\bibitem[PP16]{Pisante2016}
	Adriano Pisante and Fabio Punzo.
	\newblock Allen-{C}ahn approximation of mean curvature flow in {R}iemannian
	manifolds {I}, uniform estimates.
	\newblock {\em Ann. Sc. Norm. Super. Pisa Cl. Sci. (5)}, 15:309--341, 2016.
	
	\bibitem[PR03]{PacardRitore}
	Frank Pacard and Manuel Ritor\'e.
	\newblock From constant mean curvature hypersurfaces to the gradient theory of
	phase transitions.
	\newblock {\em J. Differential Geom.}, 64(3):359--423, 2003.
	
	\bibitem[Pre09]{Presutti2009}
	Errico Presutti.
	\newblock {\em Scaling limits in statistical mechanics and microstructures in
		continuum mechanics}.
	\newblock Theoretical and Mathematical Physics. Springer, Berlin, 2009.
	
	\bibitem[PX09]{PacardXu}
	F.~Pacard and X.~Xu.
	\newblock Constant mean curvature spheres in {R}iemannian manifolds.
	\newblock {\em Manuscripta Math.}, 128(3):275--295, 2009.
	
	\bibitem[vdW88]{vanderWaals}
	J.~D. van~der Waals.
	\newblock {\em On the continuity of the gaseous and liquid states}.
	\newblock Studies in Statistical Mechanics, XIV. North-Holland Publishing Co.,
	Amsterdam, 1988.
	\newblock Translated from the Dutch, Edited and th an introduction by J. S.
	Rowlinson.
	
	\bibitem[ZZ19a]{arXiv:1808.03527v1}
	Xin Zhou and Jonathan~J. Zhu.
	\newblock Existence of hypersurfaces with prescribed mean curvature i - generic
	min-max.
	\newblock {\em arXiv:1808.03527 [math.DG]}, 2019.
	
	\bibitem[ZZ19b]{MR4011704}
	Xin Zhou and Jonathan~J. Zhu.
	\newblock Min-max theory for constant mean curvature hypersurfaces.
	\newblock {\em Invent. Math.}, 218(2):441--490, 2019.
	
\end{thebibliography}

\end{document}